\documentclass[aos,MSNbibl,nameyear,dvips]{arximspdf}

\doi{10.1214/13-AOS1170} 
\volume{41}
\issue{6}
\pubyear{2013}
\firstpage{2852}
\lastpage{2876}

\makeatletter
\newcommand{\rrvert}{\vert}
\newcommand{\llvert}{\vert}
\newtheorem{theorem}{Theorem}
\newproclaim{condition}{Condition}
\newtheorem{corollary}{Corollary}
\newtheorem{lemma}{Lemma}
\makeatother

\begin{document}
\begin{frontmatter}

\title{Confidence sets in sparse regression}
\runtitle{Sparse regression}
\pdftitle{Confidence sets in sparse regression}

\begin{aug}
\author[A]{\fnms{Richard} \snm{Nickl}\corref{}\ead[label=e2]{r.nickl@statslab.cam.ac.uk}}
\and
\author[B]{\fnms{Sara} \snm{van de Geer}\ead[label=e1]{geer@stat.math.ethz.ch}}
\runauthor{R. Nickl and S. van de Geer}
\affiliation{University of Cambridge and ETH Z\"urich}
\dedicated{Dedicated to the memory of Yuri I. Ingster}
\address[A]{Statistical Laboratory\\
Department of Pure Mathematics\\
\quad and Mathematical Statistics\\
University of Cambridge\\
CB3 0WB Cambridge\\
United Kingdom\\
\printead{e2}} 
\address[B]{Seminar for Statistics\\
ETH Z\"urich\\
R\"amistrasse 101\\
8092 Z\"urich\\
Switzerland\\
\printead{e1}}
\end{aug}

\received{\smonth{2} \syear{2013}}
\revised{\smonth{7} \syear{2013}}

%
\begin{abstract}
The problem of constructing confidence sets in the high-dimen\-sional
linear model with $n$ response variables and $p$ parameters, possibly
$p \ge n$, is considered. Full honest adaptive inference is possible if
the rate of sparse estimation does not exceed $n^{-1/4}$, otherwise
sparse adaptive confidence sets exist only over strict subsets of the
parameter spaces for which sparse estimators exist. Necessary and
sufficient conditions for the existence of confidence sets that adapt
to a fixed sparsity level of the parameter vector are given in terms
of minimal $\ell^2$-separation conditions on the parameter space. The
design conditions cover common coherence assumptions used in models for
sparsity, including (possibly correlated) sub-Gaussian designs.
\end{abstract}

%
\begin{keyword}[class=AMS]
\kwd[Primary ]{62J05}
\kwd[; secondary ]{62G15}
\end{keyword}
\begin{keyword}
\kwd{Composite testing problem}
\kwd{high-dimensional inference}
\kwd{detection boundary}
\kwd{quadratic risk estimation}
\end{keyword}

\end{frontmatter}

\section{Introduction}\label{sec1}
 Consider the linear model
%
\begin{equation}
\label{model} Y = X \theta+ \varepsilon,
\end{equation}
where $X$ is a $n \times p$ matrix, $\theta\in\mathbb R^p$, potentially
$p>n$, and where $\varepsilon$ is a $n\times1$ vector consisting of
i.i.d. Gaussian noise independent of $X$, with mean zero and known
variance standardised to one. To develop the main ideas, let us assume
for the moment that the matrix $X$ consists of i.i.d. $N(0,1)$ Gaussian
entries $(X_{ij})$, reflecting a prototypical high-dimensional model,
such as those encountered in compressive sensing; our main results hold
for more general design assumptions that we introduce and discuss in
detail below.

We denote by $P_\theta$ the law of $(Y, X)$, by $E_\theta$ the
corresponding expectation, and will omit the subscript $\theta$ when no
confusion may arise. For the asymptotic analysis we shall let $\min
(n,p)$ tend towards infinity, and the $o,O$-notation is to be
understood accordingly. Let $B_0(k)$ be the $\ell^0$-``ball'' of radius
$k$ in $\mathbb R^p$, so all vectors in $\mathbb R^p$ with at most $k
\le p$ nonzero entries. As common in the literature on high-dimensional
models, we shall consider $p$ potentially greater than $n$ but signals
$\theta$ that are \textit{sparse} in the sense that $\theta\in B_0(k)$
for some $k$ significantly smaller than $p$. We parameterise $k$ as
\[
k\equiv k(\beta) \sim p^{1-\beta},\qquad 0<\beta<1.
\]
The parameter $\beta$ measures the sparsity of the signal: if $\beta$
is close to one, only very few of the $p$ coefficients of $\theta$ are
nonzero. If $\beta\in(0,1/2]$, one speaks of the moderately sparse case
and for $\beta\in(1/2,1]$ of the highly sparse case. We include the
case $\beta=1$ where, by convention, $k \equiv \operatorname{const}
\times p^{0} = \operatorname{const}$.

A sparse adaptive estimator $\hat\theta\equiv\hat\theta_{np} = \hat
\theta(Y,X)$ for $\theta$ achieves for every $n$, every $k \le p$, some
universal constant $c$ and with high $P_\theta$-probability, the risk bound
%
\begin{equation}
\label{spast} \|\hat\theta- \theta\|^2 \le c \log p \times
\frac{k}{n},
\end{equation}
uniformly for all $\theta\in B_0(k)$. Here $\|\cdot\| \equiv\|\cdot\|
_2$ denotes the standard Euclidean norm on $\mathbb R^p$, with inner
product $\langle\cdot, \cdot\rangle$. Such estimators exist (see
Corollary~\ref{spadest} below, for example)---they attain the risk of an
estimator that would know the positions of the $k$ nonzero
coefficients, with the typically mild penalty of $\log p$. The
literature on such estimators is abundant; see, for instance,
\citet{CT07}, \citet{BRT09}, and the monograph
\citet{BvdG2011}, where many further references can be found.

We are interested in the question of whether one can construct a
confidence set for $\theta$ that takes inferential advantage of
sparsity as in (\ref{spast}). Most of what follows applies as well to
the related problem of constructing confidence sets for $X\theta$---we
discuss this briefly at the end of the \hyperref[sec1]{Introduction}. A confidence set
$C\equiv C_{np}$ is a~random subset of $\mathbb R^p$---depending only
on the sample $Y$, $X$ and on a significance level $0<\alpha<1$---that
we require to contain the true parameter $\theta$ with at least a
prescribed probability $1-\alpha$. Our positive results rely on the
in many ways natural universal assumption $\theta\in B_0(k_1)$, with $k_1$
a minimal sparsity degree such that consistent estimation is possible.
Formally,
\[
k_1 \sim p^{1-\beta_1},\qquad \beta_1 \in(0,1);\qquad
k_1 = o(n/\log p),
\]
so that the risk bound in (\ref{spast}) converges to zero for $k=k_1$.
Our statistical procedure should have coverage over signals that are at
least $k_1$-sparse. Given $0<\alpha<1$, a~level $1-\alpha$ confidence set
$C$ should then be \textit{honest} over $B_0(k_1)$,
%
\begin{equation}
\label{honest} \liminf_{\min(n,p) \to\infty}\ \inf_{\theta\in B_0(k_1)}
P_{\theta}(\theta\in C) \ge1- \alpha.
\end{equation}
Moreover, the Euclidean diameter $|C|_2$ of $C$ should satisfy that for
every $\alpha'>0$ there exists a universal constant $L$ such that for
every $0<k \le k_1$,
%
\begin{equation}
\label{adapt} \limsup_{\min(n,p) \to\infty}\ \sup_{\theta\in B_0(k)}
P_\theta \biggl(|C|^2_2 > L \log p \times
\frac{k}{n} \biggr) \le\alpha'.
\end{equation}
Such a confidence set would cover the true $\theta$ with prescribed
probability and would shrink at an optimal rate for $k$-sparse signals
without requiring knowledge of the position of the $k$ nonzero coefficients.

A first attempt to construct such a confidence set, inspired by
\citet{L89}, \citet{BD98}, \citet{B04} in nonparametric
regression problems, is based on estimating the accuracy of estimation
in (\ref{spast}) directly via sample splitting. Heuristically the idea
is to compute a sparse estimator $\tilde\theta$ based on the first
subsample of $(Y,X)$ and to construct a confidence set centred at
$\tilde\theta$ based on the risk estimate
\[
\frac{1}{n}(Y-X \tilde\theta){}^T(Y-X\tilde\theta)-1
\]
based on $Y,X$ from the other subsample.

\begin{theorem} \label{n4}
Consider the model (\ref{model}) with i.i.d. Gaussian design $X_{ij}
\sim N(0,1)$ and assume $k_1=o(n/\log p)$. There exists a confidence
set $C$ that is honest over $B_0(k_1)$ in the sense of (\ref{honest})
and which satisfies, for any $k \le k_1$, and uniformly in $\theta\in
B_0(k)$,
\[
|C|_2^2 = O_P \biggl(\log p \times
\frac{k}{n} + n^{-1/2} \biggr).
\]
\end{theorem}

In fact, we prove Theorem~\ref{n4} for general correlated designs
satisfying Condition~\ref{subgauss} below. As a consequence, in such
situations full adaptive inference is possible if the rate of sparse
estimation in (\ref{spast}) is not desired to exceed~$n^{-1/4}$.

One may next look for estimates of $\|\tilde\theta-\theta\|$ that have
a better accuracy than just of order $n^{-1/4}$. In nonparametric
estimation problems this has been shown to be possible; see
\citet{HL02}, \citet{JL03}, \citet{CL06},
\citet{RV06}, \citet{BN12}. Translated to high-dimensional
linear models, the accuracy of these methods can be seen to be of order
$p^{1/4}/\sqrt n$, which for $p \ge n$ is of larger order of magnitude
than $n^{-1/4}$ and hence of limited interest.

Indeed, our results below will show that the rate $n^{-1/4}$ is
intrinsic to high-dimensional models: for $p \ge n$ a confidence set
that simultaneously satisfies (\ref{honest}) and adapts at \textit{any}
rate $\sqrt{(k \log p)/n} = o(n^{-1/4})$ in (\ref{adapt}) does not
exist. Rather one then needs to remove certain `critical regions' from
the parameter space in order to construct confidence sets. This is so
despite the existence of estimators satisfying (\ref{spast});
\textit{the construction of general sparse confidence sets is thus a
qualitatively different problem than that of sparse estimation.}

To formalise these ideas, we take the separation approach to adaptive
confidence sets introduced in \citet{GN10}, \citet{HN11},
\citet{BN12} in the framework of nonparametric function
estimation. We shall attempt to make honest inference over maximal
subsets of $B_0(k_1)$ where $k_1$ is given a priori as above, in a way
that is adaptive over the submodel of sparse vectors $\theta$ that
belong to $B_0(k_0)$,
\[
k_0 \sim p^{1-\beta_0},\qquad k_0 < k_1,\qquad \beta_0>\beta_1.
\]
By tracking constants in our proofs, we could include $\beta_0 = \beta
_1$ too if $k_0 \le ck_1$ for $c>0$ a small constant without changing
our findings. However, assuming $k_0=o(k_1)$ results in a considerably
cleaner mathematical exposition.

We shall remove those $\theta\in B_0(k_1)$ that are too close in
Euclidean distance to $B_0(k_0)$, and consider
%
\begin{equation}
\label{sepclass} \tilde B_0(k_1, \rho) = \bigl\{\theta\in
B_0(k_1)\dvtx  \bigl\|\theta- B_0(k_0)\bigr\|
\ge\rho \bigr\},
\end{equation}
where $\rho= \rho_{np}$ is a separation sequence and where $\|\theta
-Z\|= \inf_{z \in Z}\|\theta-z\|$ for any $Z \subset\mathbb R^p$. Thus,
if $\theta\notin B_0(k_0)$, we remove the $k_0$ coefficients $\theta_j$
with largest modulus $|\theta_j|$ from $\theta$, and require a lower
bound on the $\ell^2$-norm of the remaining subvector. In other words,
if $|\theta_{(1)}| \le\cdots\le|\theta_{(j)}| \le\cdots\le
|\theta_{(p)}|$ are any order statistics of $\{|\theta_j|\dvtx  j=1,
\ldots, p\}$, then
\[
\bigl\|\theta-B_0(k_0)\bigr\|^2 = \sum_{j=1}^{p-k_0} \theta_{(j)}^2
\]
needs to exceed $\rho^2$. Defining the new model
\[
\Theta(\rho) = B_0(k_0) \cup\tilde B_0(k_1,
\rho),
\]
we now require, instead of (\ref{honest}) and (\ref{adapt}), the weaker
coverage property
%
\begin{equation}
\label{h2} \liminf_{\min(n,p) \to\infty}\ \inf_{\theta\in\Theta(\rho_{np})}
P_{\theta}(\theta\in C_{np}) \ge1- \alpha
\end{equation}
for any $0<\alpha<1$, as well as, for some finite constant $L>0$,
%
\begin{equation}
\label{ad1} \limsup_{\min(n,p) \to\infty}\ \sup_{\theta\in B_0(k_0)}
P_\theta \biggl(|C_{np}|^2_2 > L \log
p \times\frac{k_0}{n} \biggr) \le\alpha'
\end{equation}
and
%
\begin{equation}
\label{ad2} \limsup_{\min(n,p) \to\infty}\ \sup_{\theta\in\tilde B_0(k_1, \rho
_{np}) }
P_\theta \biggl(|C_{np}|^2_2 > L \log
p \times\frac{k_1}{n} \biggr) \le\alpha'
\end{equation}
and search for minimal assumptions on the separation sequence $\rho
_{np}$. Note that any confidence set $C$ that satisfies (\ref{honest})
and (\ref{adapt}) also satisfies the above three conditions for any
$\rho\ge0$, so if one can prove the necessity of a lower bound on the
sequence~$\rho_{np}$, then one disproves in particular the existence of
adaptive confidence sets in the stronger sense of (\ref{honest}) and
(\ref{adapt}).

The following result describes our findings under the conditions of
Theorem~\ref{n4}, but now requiring adaptation to $B_0(k_0)$ at
estimation rate $\sqrt{(k_0 \log p)/n}$ faster than $n^{-1/4}$ or, what
is the same, assuming
\[
k_0 = o(\sqrt n / \log p).
\]
When specialising to the high-dimensional case $p \ge n$ this
automatically forces $\beta_0>1/2$. We require coverage over moderately
sparse alternatives ($\beta_1\le1/2$); the cases $\beta_1 >1/2$, $p\le
n$ as well more general design assumptions will be considered below.

\begin{theorem} \label{gauss}
Consider the model (\ref{model}) with i.i.d. Gaussian design $X_{ij}
\sim N(0,1)$ and $p \ge n$. For $0<\beta_1 \le1/2<\beta_0 \le1$ and
$k_0< k_1$ as above assume
\[
k_0 = o(\sqrt n / \log p),\qquad k_1= o (n/\log p).
\]
An honest adaptive confidence set $C_{np}$ over $\Theta(\rho_{np})$ in
the sense of (\ref{h2}), (\ref{ad1}), (\ref{ad2}) exist if and only if
$\rho_{np}$ exceeds, up to a multiplicative universal constant,
$n^{-1/4}$, which is the minimax rate of testing between the composite
hypotheses
%
\begin{equation}
\label{compo} H_0\dvtx  \theta\in B_0(k_0)\quad
\mbox{vs.}\quad H_1\dvtx  \theta\in\tilde B_0(k_1,
\rho_{np}).
\end{equation}
\end{theorem}


The question arises whether insisting on exact rate adaptation in
(\ref{ad1}) is crucial in Theorem~\ref{gauss} or whether some mild
`penalty' for adaptation (beyond $\log p$) could be paid to avoid
separation conditions ($\rho>0$). The proof of Theorem~\ref{gauss}
implies that requiring $|C|_2^2$ in (\ref{ad1}) to shrink at any rate
$o(n^{-1/2})$ that is possibly slower than $(k_0 \log p)/n$ but still
$o((k_1 \log p)/n)$ does not alter the conclusion of necessity of
separation at rate $\rho\simeq n^{-1/4}$ in Theorem~\ref{gauss}. In
particular, for $p \ge n$, Theorem~\ref{n4} cannot be improved if one
wants adaptive confidence sets that are honest over all of $B_0(k_1)$.


Theorem~\ref{gauss} and our further results below show that sparse
$o(n^{-1/4})$-adaptive confidence sets exist precisely over those
parameter subspaces of $B_0(k_1)$ for which the degree of sparsity is
asymptotically detectable. Sparse adaptive confidence sets solve the
composite testing problem (\ref{compo}) in a minimax way, either
implicitly or explicitly. Theorem~\ref{gauss} reiterates the findings
in \citet{HN11} and \citet{BN12} that adaptive confidence
sets exist over parameter spaces for which the structural property one
wishes to adapt to--in the present case, sparsity---can be detected
from the sample.

The paper \citet{ITV10}, where the testing problem (\ref{compo})
is considered with simple $H_0\dvtx  \theta=0$, is instrumental for our
lower bound results. Our upper bounds show that a minimax test for
the composite problem (\ref{compo}) exists without requiring stronger
separation conditions than those already needed in the case of
$H_0\dvtx \theta=0$, and under general correlated design assumptions.
In the setting of Theorem~\ref{gauss} the tests are based on rejecting
$H_0$ if $T_n$ defined by
%
\begin{equation}
\label{testdef} t_n\bigl(\theta'\bigr) =
\frac{1}{\sqrt{2n}} \sum_{i=1}^n \bigl[
\bigl(Y_i-\bigl(X\theta '\bigr)_i
\bigr)^2-1\bigr],\qquad T_n = \inf_{\theta' \in B_0(k_0)}
\bigl|t_n\bigl(\theta'\bigr)\bigr|
\end{equation}
exceeds a critical value. In practice, the computation of $T_n$
requires a convex relaxation of the minimisation problem as is standard
in the construction of sparse estimators. The
proofs that such minimum tests are minimax optimal are based on ratio
empirical process techniques, particularly Lemmas
\ref{AssumptionBlemma} and~\ref{peel} below, which are of independent
interest.

Our results give weakest possible conditions on the regions of the
parameter space that have to be removed from consideration in order to
obtain sparse adaptive confidence sets for $\theta$. Another separation condition
that may come to mind would be a lower bound $\gamma_{np}$ on the
smallest nonzero entry of $\theta\in B_0(k_1)$. Then
\[
\bigl\|\theta-B_0(k_0)\bigr\|^2 \ge(k_1-k_0)
\gamma^2_{np}
\]
and if one considers, for example, moderately sparse $\beta_1<1/2$, and
$p \ge n, k_0 = o(k_1)$, the lower bound required on $\gamma_{np}$ for
Theorem~\ref{gauss} to apply is in fact $o(n^{-1/2})$. A sparse
estimator will not be able to detect nonzero coefficients of such size,
rather one needs tailor-made procedures as presented here, and similar in spirit to results
in sparse signal detection [\citet{ITV10}, \citet{ACCP11}].

Our results concern confidence sets for the parameter vector $\theta$
itself in the Euclidean norm $\|\cdot\|$. Often, instead of on $\theta
$, inference on $Z\theta$ is of interest, where $Z$~is a $m \times p$
prediction vector. If
\[
\|Z\theta\| \ge c \|\theta\|\qquad \forall\theta\in B_0(k_1)
\]
with high probability, including the important case $Z=X$ under the
usual coherence assumptions on the design matrix $X$, then any honest
confidence set for $Z\theta$ can be used to solve the testing problem
(\ref{contr}) below, so that lower bounds for sparse confidence sets
for $\theta$ carry over to lower bounds for sparse confidence sets for
$Z\theta$. In contrast, for regular fixed linear functionals of $\theta
$ such as low-dimensional projections, the situation may be different:
for instance, in the recent papers of \citet{ZZ11},
\citet{vdGBR2013} and \citet{JM13} one-dimensional confidence
intervals for a fixed element $\theta_j$ in the vector $\theta$ are
constructed.

\section{Main results}

A heuristic summary of our findings for all parameters simultaneously
is as follows: if the rate of estimation in the submodel $B_0(k_0)$ of
$B_0(k_1)$ one wishes to adapt to is faster than
%
\begin{equation}
\label{unifrat} \rho\simeq\min \biggl(n^{-1/4}, \frac{p^{1/4}}{\sqrt n}, \sqrt{
\frac
{k_1 \log p}{n}} \biggr),
\end{equation}
then separation is necessary for adaptive confidence sets to exist at
precisely this rate $\rho$. For $p \ge n$ this simply reduces to
requiring that the rate of adaptive estimation in $B_0(k_0)$ beats
$n^{-1/4}$---the natural condition expected in view of Theorem~\ref{n4}, which proves existence of honest adaptive confidence sets
when the estimation rate is $O(n^{-1/4})$.

We consider the following conditions on the design matrix $X$.

\begin{condition} \label{uncor}
Consider the model (\ref{model}) with independent and identically
distributed $(X_{ij})$ satisfying $EX_{ij}=0$, $EX^2_{ij}=1\ \forall
i,j$.
\begin{longlist}[(a)]
\item[(a)] For some $h_0>0$,
\[
\max_{1 \le j \le l \le p} E\bigl(\exp(hX_{1j}X_{1l})
\bigr) = O(1)\qquad\forall |h|\le h_0.
\]

\item[(b)] $|X_{ij}| \le b$ for some $b>0$ and all $i,j$.
\end{longlist}
\end{condition}

Let next $\hat\Sigma:= X^T X / n $ denote the Gram matrix and let
$\Sigma:= E \hat\Sigma$. We will sometimes write $\| X \theta\|_n^2:=
\theta^T \hat\Sigma\theta$ to expedite notation.

\begin{condition}\label{subgauss}
In the model (\ref{model}) assume the
following:
\begin{longlist}[(a)]
\item[(a)] The matrix $X$ has independent rows, and for each $i
    \in\{1,
    \ldots, n\} $ and each $u \in\mathbb{R}^p$ with $u^T \Sigma u
    \le1$, the random variable $(Xu)_i$ is sub-Gaussian with constants
    $\sigma_0$ and $\kappa_0$:
\[
\kappa_0^2 \bigl( E \exp\bigl[ \bigl|(Xu)_i
\bigr|^2 / \kappa_0^2 \bigr] -1 \bigr) \le
\sigma_0^2\qquad \forall u^T \Sigma u \le1.
\]

\item[(b)] The smallest eigenvalue $\Lambda_{\min}^2
    \equiv\Lambda_{\min, p}^2$ of $\Sigma$ satisfies $\inf_p
    \Lambda_{\min, p}^2 >0$.
\end{longlist}
\end{condition}

Condition~\ref{uncor}(a) could be replaced by a fixed design assumption
as in Remark~4.1 in \citet{ITV10}. Condition~\ref{uncor}(b)
clearly implies Condition~\ref{uncor}(a); it also implies Condition
\ref{subgauss} with $\Sigma=I$ and universal constants $\kappa_0,
\sigma _0$: we have $(Xu)_i = \sum_{m=1}^p X_{im} u_m$ with mean zero
and independent summands bounded in absolute value by $b|u_m|$, so that
by Hoeffding's inequality $(Xu)_i$ is sub-Gaussian,
\[
P\bigl(\bigl|(Xu)_i\bigr| \ge t \bigr) \le2 e^{-t^2/2b \|u\|_2^2}
\]
and Condition~\ref{subgauss} follows, integrating tail probabilities.

\subsection{\texorpdfstring{Adaptation to sparse signals when \mbox{$p\! \ge n$}}
{Adaptation to sparse signals when p >= n}}

We first give a version of Theorem~\ref{gauss} for general (not
necessarily Gaussian) design matrices. The proofs imply that part~(B)
actually holds also for $p \le n$ and for $0<\beta_1 < \beta_0 \le1$.

\begin{theorem}[(Moderately sparse case)]\label{sp0}
Let $p \ge n$,
$0<\beta_1 \le1/2 < \beta_0 \le1$ and let $k_0 \sim p^{1-\beta_0}<k_1
\sim p^{1-\beta_1}$ such that $k_0 = o(\sqrt n / \log p)$.
\begin{longlist}[(A)]
\item[(A) Lower bound.] Assume Condition~\ref{uncor}\textup{(a)} and
    that $\log^3 p = o(n)$. Suppose for some separation sequence $\rho_{np}
    \ge0$ and some $0<\alpha$, $\alpha' <1/3$, the confidence set
    $C_{np}$ is both honest over $\Theta(\rho_{np})$ and adapts to
    sparsity in the sense of (\ref{ad1}), (\ref{ad2}). Then necessarily
\[
\liminf_{n,p} \frac{\rho_{np}}{n^{-1/4}}>0.
\]

\item[(B) Upper bound.] Assume Condition~\ref{subgauss} and $k_1 =
    o (n/\log p)$. Then for every $0<\alpha$, $\alpha'<1$ there exists a
    sequence $\rho_{np} \ge0$ satisfying
\[
\limsup_{n,p} \frac{\rho_{np}}{n^{-1/4}}<\infty
\]
and a level $\alpha$-confidence set $C_{np}$ that is honest over
$\Theta (\rho_{np})$ and that adapts to sparsity in the sense of
(\ref{ad1}), (\ref{ad2}).
\end{longlist}
\end{theorem}

We next consider restricting the maximal parameter space itself to
highly sparse $\theta\in B_0(k_1), \beta_1 > 1/2$. If the rate of
estimation in $B_0(k_1)$ accelerates beyond $n^{-1/4}$, then one can
take advantage of this fact, although separation of $B_0(k_0)$ and
$B_0(k_1)$ is still necessary to obtain sparse adaptive confidence
sets. The following result holds also for $p \le n$.

\begin{theorem}[(Highly sparse case)]\label{sp3}
Let $1/2 < \beta_1 <
\beta_0 \le1$ and let $k_0 \sim p^{1-\beta_0}<k_1 \sim p^{1-\beta_1}$
such that $k_0 = o(\sqrt n / \log p)$.
\begin{longlist}[A]
\item[(A) Lower bound.] Assume Condition~\ref{uncor}\textup{(a)} and
    that $\log ^3 p = o(n)$. Suppose for some separation sequence $\rho_{np}
    \ge0$ and some $0<\alpha, \alpha' <1/3$, the confidence set
    $C_{np}$ is both honest over $\Theta(\rho_{np})$ and adapts to
    sparsity in the sense of (\ref{ad1}), (\ref{ad2}). Then necessarily
\[
\liminf_{n,p} \frac{\rho_{np}}{\min (\sqrt{\log p \times
({k_1}/{n})}, n^{-1/4} )}>0.
\]

\item[(B) Upper bound.] Assume Condition~\ref{subgauss} and that
    $k_1= o (n/\log p)$. Then for every $0<\alpha', \alpha<1$ there
    exists a sequence $\rho_{np} \ge0$ satisfying
\[
\limsup_{n,p} \frac{\rho_{np}}{\min (\sqrt{\log p \times
({k_1}/{n})}, n^{-1/4} )}<\infty
\]
and a level $\alpha$-confidence set $C_{np}$ that is honest over
$\Theta (\rho_{np})$ and that adapts to sparsity in the sense of
(\ref{ad1}), (\ref{ad2}).
\end{longlist}
\end{theorem}

\subsection{\texorpdfstring{The case $p\le n$---approaching nonparametric models}
{The case p <= n---approaching nonparametric models}}
\label{plen}

The case of highly sparse alternatives and $p \le n$ was already
considered in Theorem~\ref{sp3}, explaining the presence of $\sqrt
{(k_1 \log p)/n}$ in (\ref{unifrat}). We thus now restrict to $0 <\beta
_1 \le1/2$ and, moreover, to highlight the main ideas, also to $\beta
_0 > 1/2$ corresponding to the most interesting highly sparse
null-models. We now require from any confidence set $C_n$ the
conditions (\ref{h2}), (\ref{ad1}), (\ref{ad2}) with the
infimum/supremum there intersected with
\[
B_r(M)= \Biggl\{\theta\in\mathbb R^p\dvtx  \|\theta
\|^r_r = \sum_{j=1}^p
|\theta_j|^r \le M^r \Biggr\}.
\]
Let us denote the new conditions by (\ref{h2}$'$), (\ref{ad1}$'$),
(\ref{ad2}$'$).

\begin{theorem} \label{sp4} Assume $p \le n$, let $0 < \beta_1 \le1/2
< \beta_0 \le1, 0<M<\infty$, and let $k_0 \sim p^{1-\beta_0}<k_1 \sim
p^{1-\beta_1}$.
\begin{longlist}[A]
\item[(A) Lower bound.] Assume Condition~\ref{uncor}\textup{(a)}, and
    suppose for some separation sequence $\rho_{np} \ge0$ and some
    $0<\alpha, \alpha'<1/3$, the confidence set $C_{np}$ is both honest
    over $\Theta (\rho_{np}) \cap B_r(M)$ and adapts to sparsity in the
    sense of (\ref{ad1}$'$), (\ref{ad2}$'$). If $r=2$ or if $r=1$ and
    $p = O(n^{2/3})$, then necessarily
\[
\liminf_{n,p} \frac{\rho_{np}}{p^{1/4}n^{-1/2}}>0.
\]

\item[(B) Lower bound.] Assume Condition~\ref{uncor}\textup{(b)} holds
    and
    either $r=1$, $k_1= o (n/\log p)$ or $r=2$, $\beta_0=1$, $k_1= o (\sqrt{n/\log p})$. Then for every $0<\alpha, \alpha'<1$ there exists a
    sequence $\rho_{np} \ge0$ satisfying
\[
\limsup_{n,p} \frac{\rho_{np}}{p^{1/4}n^{-1/2}}<\infty
\]
and a level $\alpha$-confidence set $C\equiv C(n,p, b,M)$ that is
honest over $\Theta(\rho_{np}) \cap\{\theta\dvtx  \|\theta\|_r \le M\}$
and that adapts to sparsity in the sense of (\ref{ad1}$'$),
(\ref{ad2}$'$).
\end{longlist}
\end{theorem}

The rate $\rho$ in the previous theorems is related to the results in
\citet{BN12} and approaches, for $p = \operatorname{const}$, the
parametric theory, where the separation rate equals, quite naturally,
$1/\sqrt n$. This is in line with the findings in \citet{P09},
\citet{PS11} in the $p \le n$ setting, who point out that a class
of specific but common sparse estimators cannot reliably be used for
the construction of confidence sets.

\section{Proofs}

All lower bounds are proved in Section~\ref{lbsec}. The proofs of
existence of confidence sets are given in Section~\ref{comp}. Theorem
\ref{n4} is proved at the end, after some auxiliary results that are
required throughout.

\subsection{\texorpdfstring{Proof of Theorems \protect\ref{gauss} (necessity), \protect\ref{sp0}\textup{(A)},
\protect\ref{sp3}\textup{(A)}, \protect\ref{sp4}\textup{(A)}}
{Proof of Theorems 2 (necessity), 3(A), 4(A), 5(A)}} \label{lbsec}

The necessity part of Theorem~\ref{gauss} follows from Theorem
\ref{sp0}(A) since any i.i.d. Gaussian matrix satisfies
Condition~\ref{uncor}(a), and since its assumptions imply the growth
condition $\log^3 p = o(n)$. Except for the $\ell^r$-norm restrictions
of Theorem~\ref{sp4} discussed at the end of the proof, Theorems
\ref{sp0}(A) and~\ref{sp4}(A) can be joined into a single statement
with separation sequence $\min(p^{1/4}n^{-1/2}, n^{-1/4})$, valid for
every $p$. We thus have to consider, for all values of $p$, two cases:
the moderately sparse case $\beta_1<1/2$ with separation lower bound
$\min(p^{1/4}n^{-1/2}, n^{-1/4})$, and the highly sparse case $\beta_1
> 1/2$ with separation lower bound $\min((\log p \times(k_1/n))^{1/2},
n^{-1/4})$. Depending on the case considered,\vspace*{-1pt} denote thus by
$\rho^*=\rho^*_{np}$ either $\min((\log p \times(k_1/n))^{1/2},
n^{-1/4})$ or $\min (p^{1/4} n^{-1/2}, n^{-1/4} )$.

The main idea of the proof follows the mechanism introduced in
\citet{HN11}. Suppose by way of contradiction that $C$ is a
confidence set as in the relevant theorems, for some sequence
$\rho=\rho_{np}$ such that
\[
\liminf_{n,p} \frac{\rho}{\rho^*}=0.
\]
By passing to a subsequence, we may replace the $\liminf$ by a proper
limit, and we shall in what follows only argue along this subsequence
$n_k\equiv n$. We claim that we can then find a further sequence $\bar
\rho_{np} \equiv\bar\rho, \rho^*_{np} \ge\bar\rho_{np} \ge\rho
_{np}$, s.t.
%
\begin{equation}
\label{squeeze} \sqrt{\log p \times\frac{k_0}{n}} = o(\bar\rho),\qquad \bar
\rho=o\bigl(\rho^*\bigr),
\end{equation}
that is, $\bar\rho$ can be taken to be squeezed between the rate of
adaptive estimation in the submodel $B_0(k_0)$ and the separation rate
$\rho^*$ that we want to establish as a lower bound. To check that this
is indeed possible, we need to verify that $(\log p \times
(k_0/n))^{1/2}$ is of smaller order than any of the three terms
\[
\sqrt{\log p \times\frac{k_1}{n}},\qquad  p^{1/4} n^{-1/2},\qquad n^{-1/4}
\]
appearing in $\rho^*$. This is obvious for the first in view of the
definition of $k_0, k_1$ ($\beta_1 < \beta_0$); follows for the second
from $\beta_0>1/2$; and follows for the third from our assumption $k_0
= o(\sqrt n / \log p)$ [automatically verified in Theorem~\ref{sp4}(A)
as $p \le n$, $\beta_0>1/2$].

For such a sequence $\bar\rho$ consider testing
\[
H_0\dvtx  \theta= 0\quad\mbox{vs.}\quad H_1\dvtx  \theta\in\tilde
B_0(k_1, \bar\rho).
\]
Using the confidence set $C$, we can test $H_0$ by $\Psi= 1\{C \cap H_1
\neq\varnothing\}$---we reject $H_0$ if $C$ contains any of the
alternatives. The type two errors satisfy
\[
\sup_{\theta\in H_1} E_\theta(1-\Psi) =\sup
_{\theta\in H_1} P_\theta (C \cap H_1 = \varnothing)
\le\sup_{\theta\in H_1} P_\theta(\theta \notin C) \le\alpha+ o(1)
\]
by coverage of $C$ over $H_1 \subset\Theta(\rho)$ (recall $\bar\rho
\ge\rho$). For the type one errors we have, again by coverage, since $0
\in B_0(k_0)$ for any $k_0$, using adaptivity (\ref{ad1}) and
(\ref{squeeze}), that
\[
E_0 \Psi= P_0(C \cap H_1 \neq\varnothing) \le
P_0 \bigl(0 \in C, |C|_2 \ge \bar\rho\bigr) +\alpha+ o(1) =
\alpha'+ \alpha+o(1).
\]
We conclude from $\min(\alpha', \alpha)<1/3$ that
%
\begin{equation}
\label{contr} E_0\Psi+ \sup_{\theta\in H_1}
E_\theta(1-\Psi) \le\alpha'+2 \alpha + o(1) < 1 +o(1).
\end{equation}
On the other hand, we now show
%
\begin{equation}
\label{contra} \liminf_{n,p} \inf_{\Psi}
\Bigl(E_0\Psi+ \sup_{\theta\in H_1} E_\theta (1-\Psi)
\Bigr)\ge1,
\end{equation}
a contradiction, so that
\[
\liminf_{n,p} \frac{\rho}{\rho^*}>0
\]
necessarily must be true. Our argument proceeds by deriving
(\ref{contra}) from Theorem~4.1 in\vspace*{-3pt} \citet{ITV10}. Let $0<c<1$, $b=
\frac{\bar \rho}{c \sqrt{k_1}}$$,  h= \frac{c k_1}{p}$, and note that
%
\begin{equation}
\label{relat} b^2ph =\frac{\bar\rho^2}{c} \ge\bar\rho^2,
\qquad b^2k_0 = o \bigl(b^2ph\bigr)
\end{equation}
using that $k_0 = o(k_1)$. Consider a product prior $\pi$ on $\theta$
with marginal coefficients $\theta_j = b \varepsilon_j$, $j=1, \ldots,
p$, where the $\varepsilon_j$ are i.i.d. with $P(\varepsilon_j=0)=1-h,
P(\varepsilon_j = 1)=P(\varepsilon_j=-1)=h/2$. We show that this prior
asymptotically concentrates on our alternative space $H_1=\tilde
B_0(k_1, \bar\rho)$. Let $Z_j = \varepsilon^2_{j}$ and denote by
$Z_{(j)}$ the corresponding order statistics (counting ties in any
order, for instance, ranking numerically by dimension), then for any
$\delta>0$ and $n$ large enough, using (\ref{relat}),
\begin{eqnarray*}
\pi\bigl(\bigl\|\theta-B_0(k_0)\bigr\|^2 < (1+\delta)
\bar\rho^2\bigr) &=& P \Biggl(b^2\sum
_{j=1}^{p-k_0} Z_{(j)} < (1+\delta)\bar
\rho^2 \Biggr)
\\
&\le& P \Biggl(b^2\sum_{j=1}^{p}
Z_{(j)} < (1+\delta)\bar\rho^2 - b^2k_0
\Biggr)
\\
&\le& P \Biggl(b^2\sum_{j=1}^{p}
\varepsilon_j^2< \bar\rho^2 \Biggr) = \pi
\bigl(\|\theta\|^2 < \bar\rho^2\bigr),
\end{eqnarray*}
which by the proof of Lemma 5.1 in \citet{ITV10} converges to $0$
as $\min(n,p) \to\infty$. Moreover, that lemma also contains the proof
that $\pi(\theta\in B_0(k_1)) \to1$ (identifying $k$ there with our
$k_1$), which thus implies $\pi(\tilde B_0(k_1, \bar\rho)) \to1$ as
$\min(n,p) \to\infty$. The testing lower bound based on this prior,
derived in Theorem 4.1 in \citet{ITV10} (cf. particularly page
1487), then implies (\ref{contra}), which is the desired contradiction.
Finally, for Theorem~\ref{sp4}, note that the above implies immediately
that $\theta\sim\pi$ asymptotically concentrates on any fixed $\ell
^2$-ball. Moreover, $E_\pi\|\theta\|_1 = bph = o(1)$ under the
hypotheses of Theorem~\ref{sp4} when $p = O(n^{2/3})$, and likewise
$\operatorname{Var}_\pi(\|\theta\|_1) = b^2ph$, so we conclude as in
the proof of Lemma 5.1 in \citet{ITV10} that the prior
asymptotically concentrates on any fixed $\ell^1$-ball in this
situation.

\subsection{\texorpdfstring{Proofs of upper bounds: Theorems \protect\ref{gauss}
(sufficiency), \protect\ref{sp0}\textup{(B)}, \protect\ref{sp3}\textup{(B)},
\protect\ref{sp4}\textup{(B)}}
{Proofs of upper bounds: Theorems 2 (sufficiency), 3(B), 4(B), 5(B)}} \label{comp}

We first note that sufficiency in Theorem~\ref{gauss} follows from
Theorem~\ref{sp0}(B) as i.i.d. Gaussian design satisfies Condition
\ref{subgauss}. The main idea, which is the same for all theorems,
follows \citet{HN11}, \citet{BN12} to solve the composite
testing problem
%
\begin{equation}
\label{ctest} H_0\dvtx  \theta\in B_0(k_0)\quad
\mbox{vs.}\quad H_1\dvtx  \theta\in\tilde B_0(k_1,
\rho)
\end{equation}
under the pa
rameter constellations of $k_0, k_1, \rho, p, n$ relevant
in Theorems~\ref{sp0}(B), \ref{sp3}(B), \ref{sp4}(B) [and in the last
case with both hypotheses intersected with $B_r(M)$, suppressed in the
notation in what follows]. Once a minimax test $\Psi$ is available for
which type one and type two errors
%
\begin{equation}
\label{mtest} \sup_{\theta\in H_0} E_\theta\Psi_n
+ \sup_{\theta\in H_1} E_\theta (1-\Psi_n) \le\gamma
\end{equation}
can be controlled, for $n$ large enough, at any level $\gamma>0$, one
takes $\tilde\theta$ to be the estimator from (\ref{spest}) below with
$\lambda$ chosen as in Lemma~\ref{AssumptionDalemma}, and constructs
the confidence set
\begin{eqnarray*}
C_n=\cases{ \displaystyle\biggl\{\theta\dvtx  \|\theta-\tilde\theta\|_2
\le L'\sqrt{\log p \frac{k_0}{n}} \biggr\},&\quad if $
\Psi_n=0$,
\vspace*{4pt}\cr
\displaystyle\biggl\{\theta\dvtx  \|\theta-\tilde\theta\|_2
\le L'\sqrt{\log p \frac{k_1}{n}} \biggr\},&\quad if $
\Psi_n=1$.}
\end{eqnarray*}
Assuming (\ref{mtest}), we now prove that $C_n$ is honest for $B_0(k_0)
\cup\tilde B_0(k_1, \rho_{np})$ if we choose the constant $L'$ large
enough: for $\theta\in B_0(k_0)$ we have from Corollary~\ref{spadest}
below, for $L'$ large,
\begin{eqnarray*}
\inf_{\theta\in B_0(k_0)} P_\theta \{\theta\in C_n \}&
\ge& 1- \sup_{\theta\in B_0(k_0)} P_\theta \biggl\{ \|\tilde\theta-
\theta\| _2 > L'\sqrt{\log p\frac{k_0}{n}} \biggr\}
\to1
\end{eqnarray*}
as $n \to\infty$. When $\theta\in\tilde B_0(k_1, \rho_{np})$, we
have that $P_\theta \{\theta\in C_n \}$ exceeds
\[
1 - \sup_{\theta\in B_0(k_1)}P_\theta \biggl\{ \|\tilde\theta-\theta
\| _2 > L'\sqrt{\log p \frac{k_1}{n}} \biggr\} -
\sup_{\theta\in\tilde
B_0(k_1, \rho_{np})}P_\theta\{\Psi_n =0\}.
\]
The first subtracted term converges to zero for $L'$ large enough, as
before. The second subtracted term can be made less than $\gamma=\alpha
$, using (\ref{mtest}). This proves that $C_n$ is honest. We now turn
to sparse adaptivity of $C_n$: by the definition of $C_n$ we always
have $|C_n| \le L'\sqrt{\log p \times k_1/n} $, so the case $\theta
\in\tilde B_0(k_1, \rho_{np})$ is proved. If $\theta\in B_0(k_0)$,
then
\[
P_\theta \biggl\{|C_n| > L'\sqrt{\log p
\frac{k_0}{n}} \biggr\} =P_\theta\{\Psi_n =1\} \le
\alpha'
\]
by the bound on the type one errors of the test, completing the
reduction of the proof to (\ref{mtest}).

\subsubsection{\texorpdfstring{Proof of Theorem \protect\ref{sp0}\textup{(B)}}
{Proof of Theorem 3(B)}}

Throughout this subsection we impose the assumptions from Theorem
\ref{sp0}---in fact, without the restriction $p \ge n$---and with $\rho
_{np} \ge L_0 n^{-1/4}$ for some $L_0$ large enough that we will choose
below. By the arguments from the previous subsection, it suffices to
solve the testing problem (\ref{mtest}) with this choice of $\rho$, for
any $\gamma>0$. Define $t_n(\theta')$, $T_n$ as in (\ref{testdef}) and
the test $\Psi_n = 1\{T_n \ge u_\gamma\}$ where $u_\gamma$ is a
suitable fixed quantile constant such that, for every $\theta\in
B_0(k_0)$, the type one error $E_\theta\Psi_n$ is bounded by
%
\begin{eqnarray}
\label{tone} \qquad P_\theta(T_n \ge u_\gamma) &\le&
P_\theta\bigl(\bigl|t_n(\theta)\bigr| \ge u_\gamma\bigr) =
P_\theta \Biggl(\frac{1}{\sqrt{2n}} \sum_{i=1}^n
\bigl(\varepsilon_i^2-1\bigr) \ge u_\gamma
\Biggr) \le\gamma.
\end{eqnarray}
For the type two errors $\theta\in H_1$, let $\theta^*$ be a minimiser
in $T_n$ (if the infimum is not attained, the argument below requires
obvious modifications). Then
\begin{eqnarray*}
\sqrt{2n} t_n\bigl(\theta^*\bigr) &=&\sum
_{i=1}^n \bigl[\bigl(Y_i-\bigl(X
\theta^*\bigr)_i\bigr)^2-1\bigr]
\\
&=&\sum_{i=1}^n \bigl[
\bigl(Y_i-(X\theta)_i + (X\theta)_i-\bigl(X
\theta^*\bigr)_i\bigr)^2-1\bigr]
\\
&=& \sum_{i=1}^n \bigl[
\bigl(Y_i-(X\theta)_i\bigr)^2-1\bigr] + 2
\bigl\langle Y-X\theta, X\bigl(\theta -\theta^*\bigr) \bigr\rangle+ \bigl\|X\bigl(
\theta-\theta^*\bigr)\bigr\|^2,
\end{eqnarray*}
so the type two errors $E_\theta(1-\Psi_n)$ are controlled by
%
\begin{eqnarray}
\label{split} && P_\theta \Biggl(\Biggl\llvert \sum
_{i=1}^n \bigl[\bigl(Y_i-(X
\theta)_i\bigr)^2-1\bigr]+ 2\bigl\langle Y-X\theta, X
\bigl(\theta-\theta^*\bigr) \bigr\rangle\nonumber
\\[-5pt]
&&\hspace*{160pt}{}  + \bigl\|X\bigl(\theta-\theta^*\bigr)
\bigr\|^2 \Biggr\rrvert < \sqrt{2n} u_\gamma \Biggr)\nonumber
\nonumber\\[-8pt]\\[-8pt]
&&\qquad \le P_\theta \Biggl(\Biggl\llvert \sum_{i=1}^n
\bigl(\varepsilon_i^2-1\bigr) \Biggr\rrvert >
\frac{\|X(\theta-\theta^*)\|^2}{2} - \sqrt{n} u_\gamma \Biggr)
\nonumber
\\
&&\quad\qquad{} + P_\theta \biggl(\bigl\llvert 2\bigl\langle\varepsilon, X\bigl(
\theta-\theta^*\bigr) \bigr\rangle\bigr\rrvert > \frac{\|X(\theta-\theta^*)\|^2}{2} - \sqrt{n}
u_\gamma \biggr).\nonumber
\end{eqnarray}
Since $\theta^* \in B_0(k_0), \theta\in\tilde B_0(k_1, \rho)$ and
$k_0+k_1 = o(n/ \log p)$, we have, from Corollary
\ref{AssumptionAcorollary} below with $t= (k_0+k_1) \log p$ that, for
$n$ large enough and with probability at least $1-4e^{-(k_0+k_1) \log
p} \to1$,
%
\begin{eqnarray}
\label{ripl} \bigl\|X\bigl(\theta-\theta^*\bigr)\bigr\|^2 &\ge&\inf
_{\theta' \in H_0} \bigl\|X\bigl(\theta-\theta '\bigr)
\bigr\|^2 \ge c(\Lambda_{\min}) n\rho_{np}^2
\ge L' \sqrt n
\end{eqnarray}
for every $L'>0$ (choosing $L_0$ large enough). We thus restrict to
this event. The probability in the last but one line of (\ref{split})
is then bounded by
\[
P_\theta \Biggl(\Biggl\llvert \sum_{i=1}^n
\bigl(\varepsilon_i^2-1\bigr) \Biggr\rrvert > \sqrt {n}
\bigl(L'-u_\gamma\bigr) \Biggr)
\]
for $n$ large enough, which can be made as small as desired by choosing
$L' \ge4u_\gamma$, as in (\ref{tone}). Likewise, the last probability
in the display (\ref{split}) is bounded, for $n$ large enough, by
\[
P_\theta \biggl(\bigl\llvert 2\bigl\langle\varepsilon, X\bigl(\theta-
\theta^*\bigr) \bigr\rangle \bigr\rrvert > \frac{\|X(\theta-\theta^*)\|^2}{4} \biggr) \le
P_\theta \biggl(\sup_{\theta' \in H_0} \frac{2\llvert \langle\varepsilon, X(\theta
-\theta') \rrvert \rangle}{\|X(\theta-\theta')\|^2}>
\frac{1}{4} \biggr),
\]
which converges to zero for large enough separation constant $L_0$,
uniformly in $\tilde B_0(k_1, \rho)$, proved in Lemma
\ref{AssumptionBlemma} below [using the lower bound (\ref{ripl}) for
$\| X(\theta-\theta')\|$ and that $\sqrt{k_0 \log p /n}=o(n^{-1/4})$].

\subsubsection{\texorpdfstring{Proof of Theorem \protect\ref{sp3}\textup{(B)}}{Proof of Theorem 4(B)}}

Throughout this subsection we impose the assumptions from Theorem
\ref{sp3}(B), with $\rho_{np}$ exceeding $L_0 \sqrt{(k_1/n) \log p}$
for some $L_0$ large enough that we will choose below (the
$n^{-1/4}$-regime was treated already in Theorem~\ref{sp0}(B), whose
proof holds for all $p$). By the arguments from the beginning of
Section~\ref{comp}, it suffices to solve the testing problem
(\ref{mtest}) with this choice of $\rho$, for any level $\gamma>0$. Let
$\tilde\theta$ be the estimator from (\ref{spest}) below with $\lambda$
chosen as in Corollary~\ref{spadest} below, and define the test
statistic
\[
T_n = \inf_{\theta\in B_0(k_0)} \|\tilde\theta- \theta
\|^2,\qquad\Psi_n = 1 \biggl\{T_n \ge D \log p
\frac{k_1}{n} \biggr\}
\]
for $D$ to be chosen. The type one errors satisfy, uniformly in $\theta
\in H_0$, for $D$ large enough,
\[
E_\theta\Psi_n \le P_\theta \biggl(\|\tilde\theta-
\theta\|^2 \ge D\log p \frac{k_1}{n} \biggr) \to0
\]
as $\min(p,n) \to\infty$, by Corollary~\ref{spadest}. Likewise, we
bound $E_\theta(1-\Psi_n)$ under $\theta\in\tilde B_0(k_1, \rho)$, for
some $\theta^* \in B_0(k_0)$, by the triangle inequality,
\begin{eqnarray*}
P_\theta \biggl(\bigl\|\tilde\theta- \theta^*\bigr\|_2^2
< C\log p \frac{k_1}{n} \biggr) &\le& P_\theta \biggl(\|\tilde\theta-
\theta\| > \bigl\|\theta ^*-\theta\bigr\| - \sqrt{C\log p \frac{k_1}{n}} \biggr)
\\
& \le& P_\theta \biggl(\|\tilde\theta- \theta\|^2
\ge(L_0-C)\log p \frac{k_1}{n} \biggr) \to0
\end{eqnarray*}
for $L_0$ large enough, again by Corollary~\ref{spadest} below.

\subsubsection{\texorpdfstring{Proof of Theorem \protect\ref{sp4}\textup{(B)}}{Proof of Theorem 5(B)}}

Throughout this subsection we impose the assumptions from Theorem
\ref{sp4}(B), with $\rho_{np} \ge L_0 p^{1/4}/\sqrt n$ for some $L_0$
large enough that we will choose below. By the arguments from the
beginning of Section~\ref{comp}, it suffices to solve the testing
problem (\ref{mtest}) [with both hypotheses there intersected with
$B_r(M)$] for this choice of $\rho$ and any level $\gamma>0$. For
$\theta' \in \mathbb R^p$ we define the $U$-statistic
\[
U_n\bigl(\theta'\bigr) = \frac{2}{n(n-1)} \sum
_{i<k} \sum_{j=1}^p
\bigl(Y_iX_{ij}-\theta_j'\bigr)
\bigl(Y_k X_{kj}-\theta_j'
\bigr),
\]
which equals $\|n^{-1}X^TY - \theta'\|^2$ with diagonal terms ($i=k$)
removed. Then
%
\begin{eqnarray}\label{exp}
\frac{1}{n}E_\theta X^TY &=&
E_\theta \biggl(\frac{1}{n}X^TX \biggr) \theta =
\theta,
\nonumber\\[-8pt]\\[-8pt]
E_\theta Y_1X_{1j} &=&
\theta_j,\qquad E_\theta U_n\bigl(
\theta'\bigr) = \bigl\| \theta-\theta'\bigr\|^2\nonumber
\end{eqnarray}
and we define the test statistic and test as
\[
T_n = \inf_{\theta' \in B_0(k_0)} \bigl|U_n\bigl(
\theta'\bigr)\bigr|,\qquad \Psi_n = 1 \biggl\{T_n
\ge u_\gamma\frac{\sqrt p}{n} \biggr\}
\]
for $u_\gamma$ quantile constants specified below. For the type one
errors we have, uniformly in $H_0$, by Chebyshev's inequality
%
\begin{equation}
\label{varest} \qquad E_\theta\Psi_n = P_\theta
\biggl(T_n \ge u_\gamma\frac{\sqrt
p}{n} \biggr) \le
P_\theta \biggl(\bigl|U_n(\theta)\bigr| \ge u_\gamma
\frac{\sqrt
p}{n} \biggr) \le\frac{\operatorname{Var}(U_n(\theta))}{u^2_\gamma} \frac{n^2}{p}.
\end{equation}
Under $P_\theta$ the $U$-statistic $U_n(\theta)$ is fully centered [cf.
(\ref{exp})], and by standard \mbox{$U$-}statistic arguments the variance can
be bounded by $\operatorname{Var}_\theta(U_n(\theta)) \le D p/n^2$ for
some constant $D$ depending only on $M$ and\vspace*{1pt} $\max_{1\le j \le
p}EX^4_{1j} \le b^4$; see, for instance, display (6.6) in
\citet{ITV10} and the arguments preceding it. We can thus choose
$u_\gamma= u_\gamma(M, b)$ to control the type one errors in~(\ref{varest}).

We now turn to the type two errors and assume $\theta\in\tilde B_0(k_1,
\rho)$: let $\theta^*$ be a minimiser in $T_n$, then $U_n(\theta^*)$
has Hoeffding decomposition $U_n(\theta^*) = U_n(\theta ) +
2L_n(\theta^*) + \|\theta^*-\theta\|^2$ with the linear term
\[
L_n\bigl(\theta'\bigr) = \frac{1}{n} \sum
_{i=1}^n \sum_{j=1}^p
(\theta_j - Y_iX_{ij}) \bigl(
\theta_j - \theta_j'\bigr).
\]
We can thus bound the type two errors $E_\theta(1-\Psi_n)$ as follows:
\begin{eqnarray*}
P_\theta \biggl(T_n < u_\gamma\frac{\sqrt p}{n}
\biggr) & \le& P_\theta \biggl(\bigl|U_n(\theta)\bigr| +
2\bigl|L_n\bigl(\theta^*\bigr)\bigr| \ge\bigl\|\theta-\theta^*\bigr\|^2 -
u_\gamma\frac{\sqrt p}{n} \biggr)
\\
& \le& P_\theta \biggl(\bigl|U_n(\theta)\bigr| \ge\frac{\|\theta-\theta^*\|^2}{2} -
u_\gamma\frac{\sqrt p}{2n} \biggr)
\\
&&{} +P_\theta \biggl(\bigl|L_n\bigl(\theta^*\bigr)\bigr| \ge
\frac{\|\theta-\theta^*\|
^2}{4} - u_\gamma\frac{\sqrt p}{4n} \biggr).
\end{eqnarray*}
By hypothesis on $\rho_{np}$ we can find $L_0$ large enough such that
$\|\theta-\theta^*\|^2 \ge\inf_{\theta' \in H_0} \|\theta-\theta'\|^2
\ge L \sqrt p/n$ for any $L>0$, so that the first probability in the
previous display can be bounded by $P_\theta(|U_n(\theta)| > u_\gamma
\sqrt p/n)$, which involves a~fully centered $U$-statistic and can thus
be dealt with as in the case of type one errors. The critical term is
the linear term, which, by the above estimate on $\|\theta-\theta^*\|$,
is less than or equal to
\[
P_\theta \biggl(\bigl|L_n\bigl(\theta^*\bigr)\bigr| \ge
\frac{\|\theta-\theta^*\|^2}{8} \biggr) \le P_\theta \biggl(\sup_{\theta' \in H_0}
\frac{|L_n(\theta')|}{\|
\theta-\theta'\|^2} > \frac{1}{8} \biggr).
\]
The process $L_n(\theta')$ can be written as
\begin{eqnarray*}
\bigl\langle\theta- n^{-1}X^TY, \theta-
\theta' \bigr\rangle &=& \bigl\langle\theta- n^{-1}X^TX
\theta, \theta- \theta' \bigr\rangle- \bigl\langle n^{-1}
X^T\varepsilon, \theta- \theta' \bigr\rangle
\\
& =& \frac{1}{n} \bigl\langle\bigl(E_\theta X^TX -
X^TX\bigr)\theta, \theta- \theta' \bigr\rangle-
\frac{1}{n}\bigl\langle\varepsilon, X\bigl(\theta- \theta'
\bigr) \bigr\rangle
\\
& \equiv & L^{(1)}_n\bigl(\theta'\bigr) +
L_n^{(2)}\bigl(\theta'\bigr)
\end{eqnarray*}
and we can thus bound the last probability by
%
\begin{equation}
\label{2terms} P_\theta \biggl(\sup_{\theta' \in H_0}
\frac{|L^{(1)}_n(\theta')|}{\|
\theta-\theta'\|^2} > \frac{1}{16} \biggr) + P_\theta \biggl(\sup
_{\theta
' \in H_0} \frac{|L^{(2)}_n(\theta')|}{\|\theta-\theta'\|^2} > \frac
{1}{16} \biggr).
\end{equation}
To show that the probability involving the second process approaches
zero, it suffices to show that
%
\begin{equation}
P_\theta \biggl(\sup_{\theta' \in H_0} \frac{\llvert \varepsilon^TX(\theta
-\theta')/n \rrvert }{\|X(\theta-\theta')\|^2/n} >
\frac{1}{16 \Lambda
} \biggr)
\end{equation}
converges to zero, using that $\sup_{v \in B_0(k_1)}\|Xv\|_2^2/(n\|v\|
_2^2) \le\Lambda$ for some $0<\Lambda<\infty$, on events of probability
approaching one, by Lemma~\ref{AssumptionAlemma} [noting $k_0+k_1 =
o(n/\log p)$]. By Lemma~\ref{AssumptionBlemma}\vadjust{\goodbreak} this last probability
approaches zero as \mbox{$\min(n,p)\to\infty$,} for $L_0$ large enough, noting
that the lower bound on $R_t$ there is satisfied for our separation
sequence $\rho_{np}$, by Corollary~\ref{AssumptionAcorollary} and
since $(k_0/n) \log p =o(p^{1/2}/n)$ in view of $\beta_0>1/2$.
Likewise, using the preceding arguments with Lemma~\ref{peel} instead
of Lemma~\ref{AssumptionBlemma}, the probability involving the first
process also converges to zero, which completes the proof.

\subsection{Remaining proofs}

\begin{lemma} \label{AssumptionAlemma} Assume Condition~\ref{subgauss}\textup{(a)} and denote by $P$ the law of $X$. Let $\theta\in
B_0(k_1)$ and $k \in\{ 1, \ldots, p \}$. Then for some constants
$\sigma$ and $\kappa$ depending only on $\sigma_0$ and $\kappa_0$,
$C_{k,k_1,p}\equiv(k+k_1+1) \log(25p)$ and for all $t >0$,
\begin{eqnarray*}
&& P \biggl( \sup_{\theta^{\prime} \in B_0(k),   ( \theta^{\prime} -
\theta){}^T\Sigma( \theta^{\prime} - \theta) \neq0 } \biggl|{ (\theta^{\prime} - \theta){}^T \hat\Sigma(\theta^{\prime} -
\theta) \over
(\theta^{\prime} - \theta){}^T\Sigma(\theta^{\prime} - \theta)
} -1 \biggr|
\\
&&\hspace*{75pt}{} \ge
4 \sigma\sqrt{ t + C_{k,k_1,p} \over n }+ 4 \kappa{ t +
C_{k,k_1,p} \over n} \biggr) \le4
\exp[-t].
\end{eqnarray*}
\end{lemma}

\begin{corollary}\label{AssumptionAcorollary}
Let $X$ satisfy Conditions~\ref{subgauss}\textup{(a)} and~\ref{subgauss}\textup{(b)}. Let
$\sigma$, $\kappa$, $\theta$, $k$, $k_1$, $C_{k,k_1,p}$ be defined as in Lemma
\ref{AssumptionAlemma}. Suppose that $k$, $k_1$ and $t>0$ are such that
\[
\biggl( { 8C_{k,k_1,p}\over n } \vee{8 t \over n } \biggr) \le \biggl(
{1 \over4 ( \sigma\vee\kappa) } \wedge1 \biggr).
\]
Then for all $\theta\in B_0(k_1)$,
\[
P_{\theta} \biggl( \bigl(\theta^{\prime} - \theta
\bigr){}^T \hat\Sigma\bigl(\theta^{\prime} - \theta\bigr) \ge
{1 \over2} \bigl\| \theta^{\prime} - \theta\bigr\|^2
\Lambda_{\min}^2\ \forall \theta^{\prime} \in
B_0(k) \biggr) \ge1- 4 \exp[-t].
\]
\end{corollary}

\begin{pf*}{Proof of Lemma~\ref{AssumptionAlemma}}
The vector $\theta
'-\theta$ has at most $k+k_1$ nonzero entries; in the lemma we may thus
replace $\theta'-\theta$ by a fixed vector in $B_0(k+k_1)$ and take the
supremum over all $k+k_1$-sparse nonzero vectors. In abuse of notation
let us still write $\theta'$ for any such vector, and fix a set $S
\subset\{ 1, \ldots, p \}$ with cardinality $|S|=k+k_1$. Let
$\mathbb{R}_S^p:= \{ \theta\in\mathbb{R}^p\dvtx  \theta_j = 0\  \forall
j \notin S \} $. We will show, for $\bar C(t,n) \equiv(t + 2(k+k_1)
\log5)/ n$, that
\[
P \biggl( \sup_{\theta^{\prime} \in\mathbb{R}_S^p,  ( \theta^{\prime}){}^T
\Sigma\theta^{\prime} \neq0 } \biggl| { ( \theta^{\prime} ){}^T \hat\Sigma\theta^{\prime} \over
( \theta^{\prime} ){}^T \Sigma\theta^{\prime}} -1 \biggr| \ge 4
\sigma\sqrt{\bar C(t,n) } + 4\kappa{\bar C(t,n)} \biggr) \le4 \exp[-t].
\]
Since there are ${ p \choose(k+k_1)} \le p^{(k+k_1)} $ sets $S$ of
cardinality $k+k_1$, the result then follows from the union bound. To
establish the inequality in the last display, it suffices to show
%
\begin{equation}
\label{toshowequation} P \Bigl( \sup_{\theta^{\prime} \in\mathcal{B}_S } \bigl| \bigl(
\theta^{\prime} \bigr){}^T \Phi\theta^{\prime} \bigr| \ge 4
\sigma\sqrt{\bar C(t,n) } + 4 \kappa{\bar C(t,n)} \Bigr) \le4 e^{-t},
\end{equation}
where $\mathcal{B}_S:= \{( \theta^{\prime} \in\mathbb{R}_S^p\dvtx
(\theta^{\prime}){}^T \Sigma\theta^{\prime} \le1 \}$ and $\Phi:= \hat
\Sigma- \Sigma$.\vadjust{\goodbreak}

We use\vspace*{-1pt} the notation $\| X u \|_{\Sigma} ^2:= u^T \Sigma u $, $u \in
\mathbb{R}^p$, and we let for $0< \delta< 1 $, $\{ X \theta_S^l
\}_{l=1}^{N(\delta)}$ be a minimal $\delta$-covering of $(\{ X
\theta^{\prime}\dvtx   \theta^{\prime} \in\mathcal{B}_S \}, \| \cdot
\|_{\Sigma} )$. Thus, for every $\theta^{\prime} \in\mathcal{B}_S$
there is a $\theta^l = \theta_S^l (\theta^{\prime}) $ such that $ \| X
( \theta^{\prime} - \theta^l ) \|_{\Sigma} \le\delta$. Note that $\{
\theta_S^l \} \subset\mathbb{R}_S^p$. Following an idea of
\citet{loh2012}, we then have
\[
\sup_{\theta^{\prime} \in\mathcal{B}_S } \bigl| \bigl( \theta^{\prime} -
\theta_S^l \bigl(\theta^{\prime}\bigr)
\bigr){}^T \Phi\bigl( \theta^{\prime} - \theta_S^l
\bigl(\theta^{\prime} \bigr) \bigr) \bigr| \le\delta^2 \sup
_{\vartheta\in\mathcal{B}_S } \vartheta^T \Phi\vartheta
\]
and also that $ \sup_{\theta^{\prime}
\in\mathcal{B}_S } | ( \theta^{\prime} - \theta_S^l (\theta^{\prime} )
){}^T \Phi\theta | \le\delta\sup_{\vartheta\in\mathcal{B}_S } |
\vartheta^T \Phi \vartheta| $. This implies with $\delta= 1/3$ that
%
\[
\sup_{\theta' \in \mathcal B_S} \bigl|\bigl(\theta'\bigr)\Phi \theta'\bigr| \le (9/2) \max_{l=1, \dots, N(1/3)} \bigl|\bigl(\theta^l_S\bigr)\Phi \bigl(\theta_S^l\bigr)\bigr|.
\]

Condition~\ref{subgauss}(a) ensures that for some constants $\sigma$
and $\kappa$ depending only on $\sigma_0$ and $\kappa_0$, for any $u$
with $\| X u \| _\Sigma\le1 $, and
any $t >0$, it holds that
\[
P \biggl( \bigl| u^T \Phi u \bigr| \ge\sigma\sqrt{t \over n} +
\kappa{t
\over n} \biggr) \le2 \exp[-t].
\]
This follows from the fact that the $((Xu)_i)$ are
sub-Gaussian, hence, the squares $((Xu)_i^2)$ are
sub-exponential. Bernstein's inequality can therefore be used [e.g.,
\citet{BvdG2011}, Lemma 14.9]. Finally, the covering number of a
ball in $k+k_1$-dimensional space is well known. Apply, for example,
Lemma 14.27 in \citet{BvdG2011}:  $N(\delta) \le( ( 2+ \delta)/
\delta)^{k+k_1} $. If we take $\delta=1/3$, this gives $N(1/3)
\le 9^{k+k_1}$. The union bound then proves (\ref{toshowequation}).
\end{pf*}

\subsubsection{\texorpdfstring{A ratio-bound for $\theta' \mapsto\varepsilon^TX(\theta -\theta')$}
{A ratio-bound for theta' -> epsilon TX(theta - theta')}}

\begin{lemma} \label{AssumptionBlemma} Suppose that $\varepsilon\sim
N (0, I)$ is independent of $X$. Let $\delta>0$. Then for any $t \ge
\max(1/\delta, 1)$, and for $R_t = t C_0 \sqrt{k_0 \log p / n} $ where
$C_0$ is a universal constant, we have for some universal constants
$C_1$ and $C_2$,
\begin{eqnarray*}
&& P \biggl( \sup_{\theta^{\prime} \in B_0 (k_0),  \| X
(\theta- \theta ^{\prime} ) \|_n > R_t} { |\varepsilon^T X( \theta-
\theta^{\prime} )| / n \over\| X ( \theta- \theta^{\prime} )\|_n^2 }
\ge\delta \bigg| X \biggr)
\\
&&\qquad \le C_1 \exp \biggl[- { t^2 \delta^2 k_0 \log p \over C_2 } \biggr].
\end{eqnarray*}
\end{lemma}

\begin{pf}
Let $\mathcal{G}_R (\theta):= \{ \theta^{\prime}\dvtx \| X ( \theta-
\theta^{\prime} )\|_n \le R, \theta^{\prime} \in B_0 (k_0) \} $. Then,
using the bound $\log{ p \choose k_0 } \le k_0 \log p $ and, for
example, Lemma 14.27 in \citet{BvdG2011}, we have, for $H(u, B,
\|\cdot\|) = \log N(u, B, \|\cdot\|)$ the logarithm of the usual
$u$-covering number of a subset $B$ of a normed space
\begin{eqnarray}
H\bigl( u, \bigl\{ X \bigl( \theta- \theta^{\prime} \bigr)\dvtx
\theta^{\prime} \in\mathcal {G}_R (\theta) \bigr\}, \| \cdot
\|_n \bigr) \le(k_0+1) \log \biggl( { 2R +u
\over u}
\biggr) + k_0 \log p,\nonumber
\\[-4pt]
\eqntext{u > 0.}
\end{eqnarray}
Indeed, if we fix the locations of the zeros, say, $\theta^{\prime} \in
B_0^{\prime} (k_0):= \{ \vartheta\dvtx   \vartheta_j= 0\  \forall j >
k_0 \} $, then $\{ X \theta^{\prime}\dvtx   \theta^{\prime} \in
B_0^{\prime } (k_0) \}$ is a $k_0$-dimensional linear space, so
\[
H\bigl( u, \bigl\{ X \theta^{\prime}\dvtx  \theta^{\prime} \in
B_0^{\prime} (k_0), \bigl\| X \theta^{\prime}
\bigr\|_n \le R \bigr\}, \| \cdot\|_n \bigr) \le
k_0 \log \biggl( { 2R + u \over u} \biggr),\qquad u > 0.
\]
Furthermore, the vector $X \theta$ is fixed,
so that $\mathcal{G}_R (\theta) $ is a subset of a ball with radius $R$
in the $(k_0 +1)$-dimensional linear space spanned by $\{ X_j
\}_{j=1}^{k_0}, X \theta$.

By Dudley's bound [see \citet{dudley1967sizes} or more recent
references such as \citet{VW96}, \citet{vandeGeer00}],
applied to the (conditional on $X$) Gaussian process $\theta'
\mapsto\varepsilon ^TX(\theta-\theta')$, and using $\int_0^c
\sqrt{\log(c/x) } \,dx = c \int_0^1 \sqrt{\log(1/x) } \,dx = c A $,
where $A$ is the constant $A=\int_0^1 \sqrt{\log(1/x) } \,dx$, we
obtain
\begin{eqnarray*}
E \Bigl[ \sup_{\theta^{\prime} \in\mathcal{G}_R (\theta) }
\bigl|\varepsilon^T X\bigl( \theta- \theta^{\prime} \bigr)\bigr| \vert X \Bigr]
&\le& C^{\prime} \int _0^R \sqrt{ nH\bigl( u, \mathcal{G}_R (\theta), \|
\cdot \|_n \bigr)} \,du
\\
&\le& C \sqrt{ 2 k_0 \log p} \sqrt n R
\end{eqnarray*}
for some universal constants $C \ge1$ and $C^{\prime}$. By the
Borell--Sudakov--Cirelson Gaussian concentration inequality [e.g.,
\citet{BLM13}], we therefore have for all $u>0$,
\[
P \biggl( \sup_{\theta^{\prime} \in\mathcal{G}_R (\theta) } \bigl|\varepsilon^T X\bigl(
\theta- \theta^{\prime} \bigr)\bigr|/n \ge C R \sqrt{2k_0
\log p \over n } + R
\sqrt{2u \over n} \bigg| X \biggr) \le\exp[-u].
\]
Substituting $u=v^2 k_0 \log p $ gives that for all $v>0$,
\[
P \biggl( \sup_{\theta^{\prime} \in\mathcal{G}_R (\theta) } \bigl|\varepsilon^T X\bigl(
\theta- \theta^{\prime} \bigr)\bigr|/n \ge(C +v) R \sqrt {2k_0 \log p \over n }
\bigg| X \biggr) \le\exp\bigl[-v^2 k_0 \log p \bigr],
\]
which implies that for all $v \ge1$,
\[
P \biggl( \sup_{\theta^{\prime} \in\mathcal{G}_R (\theta) } \bigl|\varepsilon^T X\bigl(
\theta- \theta^{\prime} \bigr)\bigr|/n \ge2v CR \sqrt{2k_0
\log p \over n }
\bigg| X \biggr) \le\exp\bigl[-{v^2 k_0 \log p } \bigr].
\]
Now insert the peeling device [see \citet{Alexander85}, the
terminology coming from \citet{vandeGeer00}, Section~5.3]. Let
$R_t:= 8 C t \sqrt{2 k_0 \log p / n } $. We then have
\begin{eqnarray*}
&& P \biggl( \sup_{\theta^{\prime} \in B_0 (k_0),  \| X (\theta-
\theta ^{\prime} ) \|_n > R_t} { |\varepsilon^T X( \theta-
\theta^{\prime} )| / n \over\| X ( \theta- \theta^{\prime}) \|_n^2 }
\ge\delta \bigg| X \biggr)
\\
&&\qquad \le\sum_{s=1}^{\infty} P \Bigl(\sup
_{\theta^{\prime} \in\mathcal
{G}_{2^s R_t} (\theta) } \bigl|\varepsilon^T X \bigl(\theta-
\theta^{\prime} \bigr)\bigr|/n \ge \delta2^{2(s-1)} R_t^2
\Big| X \Bigr)
\\
&&\qquad = \sum_{s=1}^{\infty} P \biggl(\sup
_{\theta^{\prime} \in\mathcal
{G}_{2^s R_t} (\theta) } \bigl|\varepsilon^T X \bigl(\theta-
\theta^{\prime} \bigr)\bigr|/n \ge 2^s R_t \times2C
\bigl(2^s t \delta\bigr) \sqrt{ 2 k_0 \log p \over n } \bigg| X
\biggr)
\\
&&\qquad \le\sum_{s=1}^{\infty} \exp\bigl[ -
2^{2s} t^2 \delta^2 k_0 \log p
\bigr] \le C_1 \exp \biggl[- { t^2 \delta^2 k_0 \log p \over C_2 } \biggr]
\end{eqnarray*}
for some universal constants $C_1$ and $C_2$, completing the proof.
\end{pf}

\subsubsection{\texorpdfstring{A ratio-bound for $\theta' \mapsto L_n^{(1)}(\theta')
\equiv\langle(E_\theta X^TX - X^TX)\theta, \theta- \theta' \rangle$}
{A ratio-bound for theta' -> L n (1)(theta') equivalent to <(E theta X TX - X TX)theta,
theta - theta'>}}

\begin{lemma} \label{peel}
We have, for every $\delta>0$, $R_t= tD_1 \sqrt{k_0 \log p/n}, t \ge
1$, some positive constants $D_1, D_2, D_3, D_4, D_5$ depending on
$\delta$, that
\[
\sup_{\theta\in B_r(M)} P_\theta \biggl(\sup_{\theta' \in B_0(k_0)\dvtx  \|
\theta-\theta'\|>R_t}
\frac{|L^{(1)}_n(\theta')|}{\|\theta-\theta'\|^2} > \delta \biggr) \le B(t,p,n),
\]
where $B(t,p,n)=D_2e^{-D_3t^2\delta^2k_0 \log p}$ under the assumptions
of Theorem~\ref{sp4}\textup{(B)}, \mbox{$r=1$}, and $B(t,p,n)=D_4 e^{-D_5 t
\delta\sqrt{n \log p}/k_1}$ under the assumptions of Theorem~\ref{sp4}\textup{(B)}, $r=2$.
\end{lemma}
\begin{pf}
The process in question is of the form
%
\begin{equation}
\quad L_n^{(1)}\dvtx \theta' \mapsto\frac{1}{n}
\sum_{i=1}^n \sum
_{j=1}^p (Z_{ij}-EZ_{ij}) \bigl(
\theta_j - \theta_j'\bigr),\qquad
Z_{ij} = \sum_{m=1}^p \theta
_m X_{im}X_{ij}.
\end{equation}
Since the $X_{ij}$ are uniformly bounded by $b$, we conclude that the
summands in $i$ of this process are uniformly bounded by
%
\begin{equation}
\label{env} 2b^2 \sum_{j=1}^p\bigl|
\theta_j-\theta_j'\bigr| \sum
_{m=1}^p |\theta_m|
\end{equation}
and the weak variances $n\operatorname{Var}_\theta (L_n^{(1)}(\theta')
)$ equal, for $\delta_{mj}$ the Kronecker delta,
%
\begin{eqnarray}
\label{varpeel}
\qquad && E\sum_{j, l} (Z_{ij}-EZ_{ij})
(Z_{il}-EZ_{il}) \bigl(\theta_j-\theta
_j'\bigr) \bigl(\theta_l-
\theta_l'\bigr)
\nonumber
\\
&&\qquad = E\sum_{j,l,m,m'} (X_{im}X_{ij}-
\delta_{mj}) (X_{im'}X_{il}-\delta
_{m'l}) \theta_m \theta_{m'} \bigl(
\theta_j-\theta_j'\bigr) \bigl(
\theta_l - \theta_l'\bigr)
\\
&&\qquad= \sum_{j,l,m,m'} D_{mjm'l}
\theta_m \theta_{m'} \bigl(\theta_j-\theta
_j'\bigr) \bigl(\theta_l -
\theta_l'\bigr) \le c\|\theta\|_2^2
\bigl\|\theta-\theta'\bigr\|_2^2,\nonumber
\end{eqnarray}
where we have used, by the design assumptions, that $D_{mjm'l}\le1$
whenever the indices $m,j,m',l$ match exactly to two distinct values,
$D_{mjm'l}\le EX^4_{11}$ if $m=l=j=m'$, and $D_{mjm'l}=0$ in all other
cases, as well as the Cauchy--Schwarz inequality. So $L_n^{(1)}$ is a
uniformly bounded empirical process $\{(P_n-P)(f_{\theta'})\}_{\theta'
\in H_0}$ given by
\[
\frac{1}{n}\sum_{i=1}^n
\bigl(f_{\theta'}(Z_i)-Ef_{\theta'}(Z_i)
\bigr),\qquad f_{\theta'}(Z_i) = \sum
_{j=1}^p\sum_{m=1}^p
\theta_m X_{im}X_{ij}\bigl(\theta_j-
\theta_j'\bigr)
\]
with variables $Z_i = (X_{i1}, \ldots, X_{ip}){}^T \in\mathbb R^p$.
Define $\mathcal F_s \equiv\{f=f_{\theta'}\dvtx  \theta' \in H_0,
\|\theta '-\theta\|^2 \le2^{s+1} \}$. We know $R_t < \|\theta-\theta'\|
\le \sqrt C$ so the first probability in (\ref{2terms}) can be bounded,
for $c'>0$ a small constant, by
\begin{eqnarray*}
&& P_\theta \biggl(\max_{s \in\mathbb Z\dvtx  c'R_t^2 \le2^s \le C} \sup
_{\theta' \in H_0, 2^s < \|\theta-\theta'\|^2 \le2^{s+1}} \frac
{|L^{(1)}_n(\theta')|}{\|\theta-\theta'\|^2} > \delta \biggr)
\\
&&\qquad \le\sum_{s \in\mathbb Z\dvtx  c' R_t^2 \le2^s \le C} P_\theta \Bigl(\sup
_{\theta' \in H_0, \|\theta-\theta'\|^2 \le2^{s+1}} \bigl|L^{(1)}_n\bigl(
\theta'\bigr)\bigr|> 2^s\delta \Bigr),
\end{eqnarray*}\vspace*{-5pt}
\[
\sum_{s \in\mathbb Z\dvtx  c' R_t^2 \le2^s \le C} P_\theta \bigl(\|
P_n-P\|_{\mathcal F_s}- E\|P_n-P\|_{\mathcal F_s}>
2^s \delta- E\| P_n-P\|_{\mathcal F_s} \bigr).
\]
Moreover, $\mathcal F_s$ varies in a linear space of measurable
functions of dimension $k_0$, so we have, from $\log{ p \choose k_0 }
\le k_0 \log p $ and from Theorem 2.6.7 and Lemma 2.6.15 in
\citet{VW96}, that
\[
H\bigl(u, \mathcal F_s, L^2(Q)\bigr) \lesssim
k_0 \log(AU/u)+ k_0 \log p,\qquad 0<u<UA
\]
for some universal constant $A$ and envelope bound $U$ of $\mathcal
F_s$. Using (\ref{env}), if $\theta, \theta'$ are bounded in $\ell^1$
by $M$, we can take $U$ a large enough fixed constant depending on $M,
b$ only, and if $k_0$ is constant, we can take $U=\max(k_1\sqrt
{2^s},1)$ since $\|\theta-\theta'\|_1 \le\sqrt{k_1} \|\theta-\theta'\|
_2$. A standard moment bound for empirical processes under a uniform
entropy condition [e.g., Proposition 3 in \citet{GN09a}] then
gives, using~(\ref{varpeel}),
%
\begin{equation}
\label{mombd} E\|P_n-P\|_{\mathcal F_s} \lesssim\sqrt{
\frac{2^s k_0}{n} \log p} + \frac{U k_0 \log p}{n},
\end{equation}
which is, under the maintained hypotheses, of smaller order than $2^s
\delta$ precisely for those $s$ such that $R_t^2 \simeq(k_0/n) \log p
\lesssim2^s$. The last sum of probabilities can thus be bounded, for
$D_1$ large enough and $c_0$ some positive constant, by
\[
\sum_{s \in\mathbb Z\dvtx  c'R_t^2 \le2^s \le C} P_\theta \bigl(n
\|P_n-P\| _{\mathcal F_s}- nE\|P_n-P\|_{\mathcal F_s}>
c_0n2^s \delta \bigr),
\]
to which we can apply Talagrand's inequality [\citet{T96}] [as at
the end of the proof of Proposition 1 in \citet{BN12}], to obtain
the bound
\[
\sum_{s \in\mathbb Z\dvtx  c'R_t^2 \le2^s \le C} \exp \biggl\{-\delta^2
\frac{c_0^2n^2(2^s)^2}{n2^{s+1} + n UE \|P_n-P\|_{\mathcal F_s} +
Uc_0n2^s\delta} \biggr\}.
\]
Using (\ref{mombd}), this gives the desired bound $D_2e^{-D_3t^2\delta
^2 k_0 \log p}$ when the envelope $U$ is constant, and the bound
$B(t,p,n)=D_4 e^{-D_5t \delta(n \log p)^{1/2}/k_1}$ when the envelope
is $U= \max(k_1\sqrt{2^s},1)$ (with $k_0$ constant), completing the proof.
\end{pf}

\subsubsection{Tail inequalities for sparse estimators}\label{tail}
Recall that \mbox{$S_{\vartheta}:= \{ j\dvtx   \vartheta_j \neq0 \} $}. Let
$k_{\vartheta}:= | S_{\vartheta} | $. For $\lambda>0$, take the
estimator
%
\begin{equation}
\label{spest} \tilde\theta:= \arg\min_{\vartheta} \bigl\{ \| Y - X
\vartheta\|_2^2 / n + \lambda^2
k_{\vartheta} \bigr\}.
\end{equation}

\begin{lemma} \label{AssumptionDalemma} Let $\varepsilon\sim\mathcal
{N} (0, I)$ be independent of $X$. Take $\lambda^2 = C_3 \log p / n $,
where $C_3$ is an appropriate universal constant. Let $t \ge1$ be
arbitrary and $R_t:= \sqrt{t/n}$. Then for some universal constants
$C_4$ and $C_5$,
\[
\sup_{\theta\in B_0 (k_0) } P_{\theta} \bigl( \bigl\| X ( \tilde\theta-
\theta) \bigr\|_n^2 + \lambda^2 k_{\tilde\theta} >
2 \lambda^2 k_0 + R_t^2 \vert X
\bigr) \le C_4 \exp \biggl[ - { n R_t^2 \over C_5 } \biggr].
\]
\end{lemma}

\begin{pf}
The result follows from an oracle inequality for
least squares estimators with general penalties as given in
\citet{vdG2001}. For completeness, we present a full proof. Define
\[
\tau^2 ( \vartheta; \theta):= \bigl\| X ( \vartheta- \theta)
\bigr\|_n^2 + \lambda^2 k_{\vartheta}\quad\mbox{and}\quad\mathcal{G}_R ( \theta):= \bigl\{ \vartheta\dvtx
\tau^2 ( \vartheta) \le R \bigr\}.
\]

If $\tau^2 ( \tilde\theta; \theta) \le2 \lambda^2 k_{\theta} $, we
are done.
So suppose $\tau^2 ( \tilde\theta; \theta) > 2 \lambda^2 k_{\theta} $.
We then have
$(2/n) \varepsilon^T X( \tilde\theta- \theta) \ge\tau^2 ( \tilde
\theta, \theta) - \lambda^2 k_{\theta} \ge\tau^2 ( \tilde\theta,
\theta) / 2$.
Now again apply the peeling device:
\begin{eqnarray*}
&& P \biggl( \sup_{\tau(\vartheta; \theta) > R_t } { \varepsilon^T X(
\vartheta- \theta)/n \over\tau^2 ( \vartheta, \theta) } \ge
{1 \over4} \bigg| X \biggr)
\\
&&\qquad \le\sum_{s=1}^{\infty} P \biggl( \sup
_{\vartheta\in\mathcal{G}_{2^s
R_t} (\theta) } { \varepsilon^T X( \vartheta- \theta)/n \over\tau^2
( \vartheta, \theta) } \ge {1 \over16 }
2^{2s} R_t^2 \bigg| X \biggr).
\end{eqnarray*}
But if $\vartheta\in\mathcal{G}_R (\theta)$, we know that $\|X(
\vartheta- \theta) \|_n \le R$ and that $k_{\vartheta} \le R^2 /
\lambda^2 $. Hence, as in the proof of Lemma~\ref{AssumptionBlemma},
we know that
\[
P \biggl( \sup_{\vartheta\in\mathcal{G}_R (\theta) } \varepsilon^T X( \vartheta-
\theta)/n \ge2 CR \sqrt{2R^2 \log p \over n\lambda^2 } \bigg| X \biggr) \le\exp \biggl[-
{ C^2 R^2 \log p \over
\lambda^2 } \biggr].
\]
As $\lambda= 32 C \sqrt{2 \log p / n } $, we get
\[
P \biggl( \sup_{\vartheta\in\mathcal{G}_R (\theta) } \varepsilon^T X( \vartheta-
\theta)/n \ge{ R^2 \over16} \bigg| X \biggr) \le\exp \biggl[-
{n R^2 \over2 \times(32)^2
} \biggr].
\]
We therefore have
\[
P \biggl( \sup_{\tau(\vartheta; \theta) > R_t } { \varepsilon^T X(
\vartheta- \theta) /n \over\tau^2 ( \vartheta, \theta) } \ge
{1 \over4} \bigg| X \biggr) \le\sum_{s=1}^{\infty}
\exp \biggl[- {n 2^{2s} R_t^2 \over2 \times(32)^2 } \biggr] \le C_4 \exp \biggl[ -
{ n R_t^2 \over C_5 } \biggr]
\]
for some universal constants $C_4$ and $C_5$.
\end{pf}

\begin{corollary} \label{spadest}
Assume Condition~\ref{subgauss} and let $\varepsilon\sim\mathcal{N} (0,
I)$ be independent of~$X$. Let $\tilde\theta$ be as in (\ref{spest})
with $\lambda^2 = (C_3 \log p) / n $, where $C_3$ is as in Lemma
\ref{AssumptionDalemma}, and let $k_0 = o(n/\log p)$. Then for some
universal constants $C_6, C_7, C_8, c$, every $C\ge C_6$ and every $n$
large enough,
\[
\sup_{\theta\in B_0 (k_0) } P_{\theta} \biggl( \| \tilde\theta- \theta
\|^2 > C \frac{k_0 \log p }{n} \biggr) \le C_7 \exp \biggl[
- {k_0 \log p \over C_8} \biggr].
\]
\end{corollary}
\begin{pf}
By Lemma~\ref{AssumptionDalemma} with $R_\tau$, $\tau$ equal to a
suitable constant times $k_0 \log p$, we see first $k_{\tilde\theta}
\lesssim3k_0$ on the event on which the exponential inequality holds.
Then from Corollary~\ref{AssumptionAcorollary} with $k=3k_0$, on an
event of sufficiently large probability, $\|\tilde\theta- \theta\| _2^2
\le C(\Lambda_{\min}) \|X(\tilde\theta-\theta)\|_{n}^2$ for $n$ large
enough, so that the result follows from applying Lemma
\ref{AssumptionDalemma} again [this time to $\|X(\tilde\theta-
\theta)\| _n^2$] and from combining the bounds.
\end{pf}

\subsubsection{\texorpdfstring{Proof of Theorem \protect\ref{n4} under Condition \protect\ref{subgauss}}
{Proof of Theorem 1 under Condition 2}}

For $p, n$ fixed, the random vectors $(Y_i, X_{i1}, \ldots,
X_{ip})_{i=1}^n$ are i.i.d., and if we split the $n$ points into two
subsamples, each of size of order $n$, then we have two independent
replicates $Y^{(s)}=X^{(s)}\theta+\varepsilon^{(s)}, \hat\Sigma^{(s)}=
(X^{(s)}){}^TX^{(s)}/n, s=1,2$, of the model. In abuse of notation,
denote throughout this proof by $\tilde\theta\equiv\tilde\theta ^{(1)}$
the estimator from (\ref{spest}) based on the subsample $s=1$, with
$\lambda$ chosen as in Lemma~\ref{AssumptionDalemma}, and by $(Y,X,
\varepsilon)\equiv(Y^{(2)}, X^{(2)}, \varepsilon^{(2)})$ the variables
from the second subsample. Define
\begin{eqnarray*}
\hat R_n &=& \frac{1}{n} (Y-X\tilde\theta){}^T(Y-X
\tilde \theta) - 1
\\
& =& (\theta- \tilde\theta){}^T\hat\Sigma^{(2)} (\theta-\tilde
\theta) + \frac{2}{n} \varepsilon^TX(\theta-\tilde\theta) +
\frac{1}{n} \varepsilon^T \varepsilon-1.
\end{eqnarray*}
By independence, and conditional on $(Y^{(1)}, X^{(1)})$, we have
$E^{(2)}_\theta(\varepsilon^TX(\theta-\tilde\theta))^2 = n(\tilde
\theta-\theta){}^T\Sigma(\tilde\theta-\theta)$ and so, using Markov's
inequality,
%
\begin{equation}
\label{linearterm} \frac{2}{n} \varepsilon^TX(\theta-\tilde
\theta) = O_P \biggl(\sqrt{\frac
{(\tilde\theta-\theta){}^T\Sigma(\tilde\theta-\theta)}{ n}} \biggr).
\end{equation}
By Lemma~\ref{AssumptionDalemma}, we have $\| X^{(1)} ( \tilde\theta -
\theta) \|^2_n= O_P ((k \log p) / n ) $ and $k_{\tilde\theta} = O(k_1)$
and, hence, by Lemma~\ref{AssumptionAlemma}, also $(\tilde\theta-
\theta){}^T \Sigma( \tilde \theta- \theta) =O_P  ((k \log p)/n  )=o(1)$.
Thus, the bound in (\ref{linearterm}) is $o_P(1/\sqrt n)$ uniformly in
$B_0(k_1)$, and this will be used in the following estimate. Let
$u_\alpha$ be suitable quantile constants to be chosen below. Take as
confidence set
\[
C_n = \biggl\{\theta\in\mathbb R^p\dvtx  \|\theta-\tilde
\theta\|^2 \le2 \Lambda_{\min}^{-2} \biggl( \hat
R_n + \frac{u_\alpha}{\sqrt n} \biggr) \biggr\}.
\]
Uniformly in $\theta\in B_0(k_1)$ with $k_1 = o (n/\log p)$, we have
again by Lemma~\ref{AssumptionDalemma} that $\tilde\theta\in
B_0(2k_1)$ on events of probability approaching one, so that, using
Corollary~\ref{AssumptionAcorollary} on these events,
\begin{eqnarray*}
P_\theta(\theta\notin C_n) &=& P_\theta \biggl(\|
\theta-\tilde\theta\| ^2 > 2 \Lambda_{\min}^{-2}
\biggl(\hat R_n + \frac{u_\alpha}{\sqrt
n} \biggr) \biggr)
\\
&\le& P_\theta \biggl((\theta-\tilde\theta){}^T\hat
\Sigma^{(2)}(\theta -\tilde\theta) > \hat R_n +
\frac{u_\alpha}{\sqrt n} \biggr) +o(1)
\\
&=& P_\theta \biggl(-\frac{1}{n}\varepsilon^T
\varepsilon+1 > \frac
{u_\alpha}{\sqrt n} + \frac{2}{n} \varepsilon^T
X(\theta-\tilde\theta) \biggr) + o(1)
\\
&=& P_\theta \Biggl(\frac{-1}{\sqrt n} \sum
_{i=1}^n \bigl(\varepsilon_i^2-1
\bigr) > \bigl(1+o(1)\bigr)u_\alpha \Biggr) +o(1) \le\alpha+ o(1)
\end{eqnarray*}
for a fixed constant $u_\alpha$. Moreover, from the previous arguments
and Corollary~\ref{spadest}, we see that, for $\theta\in B_0(k)$, the
diameter $\hat R_n = O_P(\|\tilde\theta-\theta\|^2+n^{-1/2})$ is of
order $O_P(\frac{k \log p}{n} + n^{-1/2})$.

\section*{Acknowledgements} We would like to thank the Editor, Associate
Editor, two referees and Sasha Tsybakov for helpful remarks on this article.
Richard Nickl is grateful to the Cafes Florianihof and Griensteidl
in Vienna where parts of this research were carried out.








\printaddresses

\begin{thebibliography}{30}

\bibitem[\protect\citeauthoryear{Alexander}{1985}]{Alexander85}
\begin{binproceedings}[mr]
\bauthor{\bsnm{Alexander},~\bfnm{Kenneth~S.}\binits{K.~S.}}
(\byear{1985}).
\btitle{Rates of growth for weighted empirical processes}.
In \bbooktitle{Proceedings of the {B}erkeley Conference in Honor of Jerzy
  {N}eyman and {J}ack {K}iefer, {V}ol. II ({B}erkeley, {C}alif., 1983)}.
\bvolume{2}
\bpages{475--493}.
\bpublisher{Wadsworth}, \blocation{Belmont, CA}.
\bid{mr={0822047}}
\bptok{imsref}%
\end{binproceedings}
\endbibitem

\bibitem[\protect\citeauthoryear{Arias-Castro, Cand{\`e}s and
  Plan}{2011}]{ACCP11}
\begin{barticle}[mr]
\bauthor{\bsnm{Arias-Castro},~\bfnm{Ery}\binits{E.}},
  \bauthor{\bsnm{Cand{\`e}s},~\bfnm{Emmanuel~J.}\binits{E.~J.}} \AND
  \bauthor{\bsnm{Plan},~\bfnm{Yaniv}\binits{Y.}}
(\byear{2011}).
\btitle{Global testing under sparse alternatives: {ANOVA}, multiple comparisons
  and the higher criticism}.
\bjournal{Ann. Statist.}
\bvolume{39}
\bpages{2533--2556}.
\bid{doi={10.1214/11-AOS910}, issn={0090-5364}, mr={2906877}}
\bptok{imsref}%
\end{barticle}
\endbibitem

\bibitem[\protect\citeauthoryear{Baraud}{2004}]{B04}
\begin{barticle}[mr]
\bauthor{\bsnm{Baraud},~\bfnm{Yannick}\binits{Y.}}
(\byear{2004}).
\btitle{Confidence balls in {G}aussian regression}.
\bjournal{Ann. Statist.}
\bvolume{32}
\bpages{528--551}.
\bid{doi={10.1214/009053604000000085}, issn={0090-5364}, mr={2060168}}
\bptok{imsref}%
\end{barticle}
\endbibitem

\bibitem[\protect\citeauthoryear{Beran and D{\"u}mbgen}{1998}]{BD98}
\begin{barticle}[mr]
\bauthor{\bsnm{Beran},~\bfnm{Rudolf}\binits{R.}} \AND
  \bauthor{\bsnm{D{\"u}mbgen},~\bfnm{Lutz}\binits{L.}}
(\byear{1998}).
\btitle{Modulation of estimators and confidence sets}.
\bjournal{Ann. Statist.}
\bvolume{26}
\bpages{1826--1856}.
\bid{doi={10.1214/aos/1024691359}, issn={0090-5364}, mr={1673280}}
\bptok{imsref}%
\end{barticle}
\endbibitem

\bibitem[\protect\citeauthoryear{Bickel, Ritov and Tsybakov}{2009}]{BRT09}
\begin{barticle}[mr]
\bauthor{\bsnm{Bickel},~\bfnm{Peter~J.}\binits{P.~J.}},
  \bauthor{\bsnm{Ritov},~\bfnm{Ya'acov}\binits{Y.}} \AND
  \bauthor{\bsnm{Tsybakov},~\bfnm{Alexandre~B.}\binits{A.~B.}}
(\byear{2009}).
\btitle{Simultaneous analysis of lasso and {D}antzig selector}.
\bjournal{Ann. Statist.}
\bvolume{37}
\bpages{1705--1732}.
\bid{doi={10.1214/08-AOS620}, issn={0090-5364}, mr={2533469}}
\bptok{imsref}%
\end{barticle}
\endbibitem

\bibitem[\protect\citeauthoryear{Boucheron, Lugosi and Massart}{2013}]{BLM13}
\begin{bbook}[auto:STB|2013/10/14|10:36:11]
\bauthor{\bsnm{Boucheron},~\bfnm{S.}\binits{S.}},
  \bauthor{\bsnm{Lugosi},~\bfnm{G.}\binits{G.}} \AND
  \bauthor{\bsnm{Massart},~\bfnm{P.}\binits{P.}}
(\byear{2013}).
\btitle{Concentration Inequalities. A~Nonasymptotic Theory of Independence}.
\bpublisher{Oxford Univ. Press}, \blocation{London}.
\bptok{imsref}%
\end{bbook}
\endbibitem

\bibitem[\protect\citeauthoryear{B{\"u}hlmann and van~de Geer}{2011}]{BvdG2011}
\begin{bbook}[mr]
\bauthor{\bsnm{B{\"u}hlmann},~\bfnm{Peter}\binits{P.}} \AND
  \bauthor{\bparticle{van~de} \bsnm{Geer},~\bfnm{Sara}\binits{S.}}
(\byear{2011}).
\btitle{Statistics for High-Dimensional Data: Methods, Theory and Applications}.
\bpublisher{Springer}, \blocation{Heidelberg}.
\bid{doi={10.1007/978-3-642-20192-9}, mr={2807761}}
\bptok{imsref}%
\end{bbook}
\endbibitem

\bibitem[\protect\citeauthoryear{Bull and Nickl}{2013}]{BN12}
\begin{barticle}[mr]
\bauthor{\bsnm{Bull},~\bfnm{Adam~D.}\binits{A.~D.}} \AND
  \bauthor{\bsnm{Nickl},~\bfnm{Richard}\binits{R.}}
(\byear{2013}).
\btitle{Adaptive confidence sets in {$L\sp 2$}}.
\bjournal{Probab. Theory Related Fields}
\bvolume{156}
\bpages{889--919}.
\bid{doi={10.1007/s00440-012-0446-z}, issn={0178-8051}, mr={3078289}}
\bptok{imsref}%
\end{barticle}
\endbibitem

\bibitem[\protect\citeauthoryear{Cai and Low}{2006}]{CL06}
\begin{barticle}[mr]
\bauthor{\bsnm{Cai},~\bfnm{T.~Tony}\binits{T.~T.}} \AND
  \bauthor{\bsnm{Low},~\bfnm{Mark~G.}\binits{M.~G.}}
(\byear{2006}).
\btitle{Adaptive confidence balls}.
\bjournal{Ann. Statist.}
\bvolume{34}
\bpages{202--228}.
\bid{doi={10.1214/009053606000000146}, issn={0090-5364}, mr={2275240}}
\bptok{imsref}%
\end{barticle}
\endbibitem

\bibitem[\protect\citeauthoryear{Cand\`es and Tao}{2007}]{CT07}
\begin{barticle}[mr]
\bauthor{\bsnm{Cand\`es},~\bfnm{Emmanuel}\binits{E.}} \AND
  \bauthor{\bsnm{Tao},~\bfnm{Terence}\binits{T.}}
(\byear{2007}).
\btitle{The {D}antzig selector: Statistical estimation when {$p$} is much
  larger than {$n$}}.
\bjournal{Ann. Statist.}
\bvolume{35}
\bpages{2313--2351}.
\bid{doi={10.1214/009053606000001523}, issn={0090-5364}, mr={2382644}}
\bptok{imsref}%
\end{barticle}
\endbibitem

\bibitem[\protect\citeauthoryear{Dudley}{1967}]{dudley1967sizes}
\begin{barticle}[mr]
\bauthor{\bsnm{Dudley},~\bfnm{R.~M.}\binits{R.~M.}}
(\byear{1967}).
\btitle{The sizes of compact subsets of {H}ilbert space and continuity of
  {G}aussian processes}.
\bjournal{J. Funct. Anal.}
\bvolume{1}
\bpages{290--330}.
\bid{mr={0220340}}
\bptok{imsref}%
\end{barticle}
\endbibitem

\bibitem[\protect\citeauthoryear{Gin{\'e} and Nickl}{2009}]{GN09a}
\begin{barticle}[mr]
\bauthor{\bsnm{Gin{\'e}},~\bfnm{Evarist}\binits{E.}} \AND
  \bauthor{\bsnm{Nickl},~\bfnm{Richard}\binits{R.}}
(\byear{2009}).
\btitle{An exponential inequality for the distribution function of the kernel
  density estimator, with applications to adaptive estimation}.
\bjournal{Probab. Theory Related Fields}
\bvolume{143}
\bpages{569--596}.
\bid{doi={10.1007/s00440-008-0137-y}, issn={0178-8051}, mr={2475673}}
\bptok{imsref}%
\end{barticle}
\endbibitem

\bibitem[\protect\citeauthoryear{Gin{\'e} and Nickl}{2010}]{GN10}
\begin{barticle}[mr]
\bauthor{\bsnm{Gin{\'e}},~\bfnm{Evarist}\binits{E.}} \AND
  \bauthor{\bsnm{Nickl},~\bfnm{Richard}\binits{R.}}
(\byear{2010}).
\btitle{Confidence bands in density estimation}.
\bjournal{Ann. Statist.}
\bvolume{38}
\bpages{1122--1170}.
\bid{doi={10.1214/09-AOS738}, issn={0090-5364}, mr={2604707}}
\bptok{imsref}%
\end{barticle}
\endbibitem

\bibitem[\protect\citeauthoryear{Hoffmann and Lepski}{2002}]{HL02}
\begin{barticle}[mr]
\bauthor{\bsnm{Hoffmann},~\bfnm{M.}\binits{M.}} \AND
  \bauthor{\bsnm{Lepski},~\bfnm{O.}\binits{O.}}
(\byear{2002}).
\btitle{Random rates in anisotropic regression}.
\bjournal{Ann. Statist.}
\bvolume{30}
\bpages{325--396}.
\bid{doi={10.1214/aos/1021379858}, issn={0090-5364}, mr={1902892}}
\bptok{imsref}%
\end{barticle}
\endbibitem

\bibitem[\protect\citeauthoryear{Hoffmann and Nickl}{2011}]{HN11}
\begin{barticle}[mr]
\bauthor{\bsnm{Hoffmann},~\bfnm{Marc}\binits{M.}} \AND
  \bauthor{\bsnm{Nickl},~\bfnm{Richard}\binits{R.}}
(\byear{2011}).
\btitle{On adaptive inference and confidence bands}.
\bjournal{Ann. Statist.}
\bvolume{39}
\bpages{2383--2409}.
\bid{doi={10.1214/11-AOS903}, issn={0090-5364}, mr={2906872}}
\bptok{imsref}%
\end{barticle}
\endbibitem

\bibitem[\protect\citeauthoryear{Ingster, Tsybakov and Verzelen}{2010}]{ITV10}
\begin{barticle}[mr]
\bauthor{\bsnm{Ingster},~\bfnm{Yuri~I.}\binits{Y.~I.}},
  \bauthor{\bsnm{Tsybakov},~\bfnm{Alexandre~B.}\binits{A.~B.}} \AND
  \bauthor{\bsnm{Verzelen},~\bfnm{Nicolas}\binits{N.}}
(\byear{2010}).
\btitle{Detection boundary in sparse regression}.
\bjournal{Electron. J. Stat.}
\bvolume{4}
\bpages{1476--1526}.
\bid{doi={10.1214/10-EJS589}, issn={1935-7524}, mr={2747131}}
\bptok{imsref}%
\end{barticle}
\endbibitem

\bibitem[\protect\citeauthoryear{Javanmard and Montanari}{2013}]{JM13}
\begin{bmisc}[auto:STB|2013/10/14|10:36:11]
\bauthor{\bsnm{Javanmard},~\bfnm{A.}\binits{A.}} \AND
  \bauthor{\bsnm{Montanari},~\bfnm{A.}\binits{A.}}
(\byear{2013}).
\bhowpublished{Confidence intervals and hypothesis testing for high-dimensional
  regression. Available at arXiv:\arxivurl{1306.3171}.}
\bptok{imsref}%
\end{bmisc}
\endbibitem

\bibitem[\protect\citeauthoryear{Juditsky and Lambert-Lacroix}{2003}]{JL03}
\begin{barticle}[mr]
\bauthor{\bsnm{Juditsky},~\bfnm{A.}\binits{A.}} \AND
  \bauthor{\bsnm{Lambert-Lacroix},~\bfnm{S.}\binits{S.}}
(\byear{2003}).
\btitle{Nonparametric confidence set estimation}.
\bjournal{Math. Methods Statist.}
\bvolume{12}
\bpages{410--428}.
\bid{issn={1066-5307}, mr={2054156}}
\bptok{imsref}%
\end{barticle}
\endbibitem

\bibitem[\protect\citeauthoryear{Li}{1989}]{L89}
\begin{barticle}[mr]
\bauthor{\bsnm{Li},~\bfnm{Ker-Chau}\binits{K.-C.}}
(\byear{1989}).
\btitle{Honest confidence regions for nonparametric regression}.
\bjournal{Ann. Statist.}
\bvolume{17}
\bpages{1001--1008}.
\bid{doi={10.1214/aos/1176347253}, issn={0090-5364}, mr={1015135}}
\bptok{imsref}%
\end{barticle}
\endbibitem

\bibitem[\protect\citeauthoryear{Loh and Wainwright}{2012}]{loh2012}
\begin{barticle}[mr]
\bauthor{\bsnm{Loh},~\bfnm{Po-Ling}\binits{P.-L.}} \AND
  \bauthor{\bsnm{Wainwright},~\bfnm{Martin~J.}\binits{M.~J.}}
(\byear{2012}).
\btitle{High-dimensional regression with noisy and missing data: Provable
  guarantees with nonconvexity}.
\bjournal{Ann. Statist.}
\bvolume{40}
\bpages{1637--1664}.
\bid{doi={10.1214/12-AOS1018}, issn={0090-5364}, mr={3015038}}
\bptok{imsref}%
\end{barticle}
\endbibitem

\bibitem[\protect\citeauthoryear{P{\"o}tscher}{2009}]{P09}
\begin{barticle}[mr]
\bauthor{\bsnm{P{\"o}tscher},~\bfnm{Benedikt~M.}\binits{B.~M.}}
(\byear{2009}).
\btitle{Confidence sets based on sparse estimators are necessarily large}.
\bjournal{Sankhy\=a}
\bvolume{71}
\bpages{1--18}.
\bid{issn={0972-7671}, mr={2579644}}
\bptok{imsref}%
\end{barticle}
\endbibitem

\bibitem[\protect\citeauthoryear{P{\"o}tscher and Schneider}{2011}]{PS11}
\begin{barticle}[mr]
\bauthor{\bsnm{P{\"o}tscher},~\bfnm{Benedikt~M.}\binits{B.~M.}} \AND
  \bauthor{\bsnm{Schneider},~\bfnm{Ulrike}\binits{U.}}
(\byear{2011}).
\btitle{Distributional results for thresholding estimators in high-dimensional
  {G}aussian regression models}.
\bjournal{Electron. J. Stat.}
\bvolume{5}
\bpages{1876--1934}.
\bid{doi={10.1214/11-EJS659}, issn={1935-7524}, mr={2970179}}
\bptok{imsref}%
\end{barticle}
\endbibitem

\bibitem[\protect\citeauthoryear{Robins and van~der Vaart}{2006}]{RV06}
\begin{barticle}[mr]
\bauthor{\bsnm{Robins},~\bfnm{James}\binits{J.}} \AND
  \bauthor{\bparticle{van~der} \bsnm{Vaart},~\bfnm{Aad}\binits{A.}}
(\byear{2006}).
\btitle{Adaptive nonparametric confidence sets}.
\bjournal{Ann. Statist.}
\bvolume{34}
\bpages{229--253}.
\bid{doi={10.1214/009053605000000877}, issn={0090-5364}, mr={2275241}}
\bptok{imsref}%
\end{barticle}
\endbibitem

\bibitem[\protect\citeauthoryear{Talagrand}{1996}]{T96}
\begin{barticle}[mr]
\bauthor{\bsnm{Talagrand},~\bfnm{Michel}\binits{M.}}
(\byear{1996}).
\btitle{New concentration inequalities in product spaces}.
\bjournal{Invent. Math.}
\bvolume{126}
\bpages{505--563}.
\bid{doi={10.1007/s002220050108}, issn={0020-9910}, mr={1419006}}
\bptok{imsref}%
\end{barticle}
\endbibitem

\bibitem[\protect\citeauthoryear{van~de Geer}{2000}]{vandeGeer00}
\begin{bbook}[auto]
\bauthor{\bparticle{van~de} \bsnm{Geer},~\bfnm{Sara~A.}\binits{S.~A.}}
(\byear{2000}).
\btitle{Empirical Processes in M-Estimation}.
\bpublisher{Cambridge Univ. Press}, \blocation{Cambridge}.
\bptok{imsref}%
\end{bbook}
\endbibitem

\bibitem[\protect\citeauthoryear{van~de Geer}{2001}]{vdG2001}
\begin{barticle}[mr]
\bauthor{\bparticle{van~de} \bsnm{Geer},~\bfnm{Sara}\binits{S.}}
(\byear{2001}).
\btitle{Least squares estimation with complexity penalties}.
\bjournal{Math. Methods Statist.}
\bvolume{10}
\bpages{355--374}.
\bid{issn={1066-5307}, mr={1867165}}
\bptok{imsref}%
\end{barticle}
\endbibitem

\bibitem[\protect\citeauthoryear{van~de Geer, B\"uhlmann and Ritov}{2013}]{vdGBR2013}
\begin{bmisc}[auto:STB|2013/10/14|10:36:11]
\bauthor{\bparticle{van~de} \bsnm{Geer},~\bfnm{Sara}\binits{S.}},
\bauthor{\bsnm{B\"uhlmann},~\bfnm{P.}\binits{P.}} \AND
\bauthor{\bsnm{Ritov},~\bfnm{Y.}\binits{Y.}}
(\byear{2013}).
\bhowpublished{On asymptotically optimal confidence regions and tests for high-dimensional
  models. Submitted. Available at arXiv:\arxivurl{1303.0518}.}
\bptok{imsref}%
\end{bmisc}
\endbibitem

\bibitem[\protect\citeauthoryear{van~der Vaart and Wellner}{1996}]{VW96}
\begin{bbook}[mr]
\bauthor{\bparticle{van~der} \bsnm{Vaart},~\bfnm{Aad~W.}\binits{A.~W.}} \AND
  \bauthor{\bsnm{Wellner},~\bfnm{Jon~A.}\binits{J.~A.}}
(\byear{1996}).
\btitle{Weak Convergence and Empirical Processes}.
\bseries{Springer Series in Statistics}.
\bpublisher{Springer}, \blocation{New York}.
\bid{mr={1385671}}
\bptok{imsref}%
\end{bbook}
\endbibitem

\bibitem[\protect\citeauthoryear{Zhang and Zhang}{2011}]{ZZ11}
\begin{bmisc}[auto:STB|2013/10/14|10:36:11]
\bauthor{\bsnm{Zhang},~\bfnm{C.~H.}\binits{C.~H.}} \AND
  \bauthor{\bsnm{Zhang},~\bfnm{S.~S.}\binits{S.~S.}}
(\byear{2011}).
\bhowpublished{Confidence intervals for low-dimensional parameters with
  high-dimensional data. 2011. Available at
  arXiv:\arxivurl{1110.2563v1}.}
\bptok{imsref}%
\end{bmisc}
\endbibitem

\end{thebibliography}
\end{document}